\documentclass[11pt,english]{article}
\usepackage{hyperref}
\hypersetup{
    colorlinks=true,
    linkcolor=blue,
    filecolor=magenta,      
    citecolor=blue,
    urlcolor=black
    }
\usepackage[sc]{mathpazo}
\usepackage[T1]{fontenc}
\usepackage{booktabs}
\usepackage{geometry}
\usepackage{color}
\usepackage{array}
\usepackage{bbding}
\usepackage{multirow}
\usepackage{algorithm}
\usepackage{algpseudocode}
\usepackage[fleqn]{amsmath}
\usepackage{exscale}
\usepackage{relsize}
\usepackage{amssymb}
\usepackage{amsthm}
\usepackage{graphicx}
\usepackage[inkscapelatex=false]{svg}
\usepackage{natbib}
\usepackage{lscape}
\usepackage{comment}
\usepackage{bm}
\usepackage[title]{appendix}
\usepackage{longtable}
\usepackage{mathtools}
\usepackage{mathrsfs}
\usepackage{chngcntr}
\usepackage{authblk}
\usepackage{tikz}
\usepackage{threeparttable}
\usepackage{makecell}
\usepackage[utf8]{inputenc}
\usepackage{subfigure}

\makeatletter
\allowdisplaybreaks

\@ifundefined{showcaptionsetup}{}{%
 \PassOptionsToPackage{caption=false}{subfig}}
\makeatother

\usepackage{babel}
\usepackage{enumitem}



\begin{document}\sloppy

\title{Electric Truck Platooning with Charging Consideration and Leader Swapping}
\author[a,b,c]{Yilang Hao}
\author[a,b]{Zhibin Chen \footnote{Corresponding author. zc23@nyu.edu.}}
\affil[a]{\small\emph{Shanghai Key Laboratory of Urban Design and Urban Science, NYU Shanghai, 567 West Yangsi Road, Pudong New District, Shanghai, 200126, China}}
\affil[b]{\small\emph{Shanghai Frontiers Science Center of Artificial Intelligence and Deep Learning, NYU Shanghai, 567 West Yangsi Road, Pudong New District, Shanghai, 200126, China}}
\affil[c]{\small\emph{Department of Civil and Urban Engineering, New York University, 6 MetroTech Center, Brooklyn, NY 11201, United States}}

\date{}
\maketitle

\begin{abstract}
Electric trucks are increasingly deployed to reduce the trucking sector’s carbon footprint. However, these vehicles often face limited driving range and range anxiety on mid to long‑haul routes, which can undermine their environmental benefits. Prior research on electric truck platooning has largely focused on a single highway corridor, failing to capture the complexity of an expansive road network across the entire country. These simplified models overlook how platooning can not only alleviate range anxiety through energy savings but also influence the selection of charging stops and even routes along a journey. Integrating platooning into electric‑truck routing over a broad network presents several challenges: 1) Range limitations and charging‑stop optimization: electric trucks must plan routes that include charging stations with varying prices and charging speed, and incorporate platoon‑induced energy savings into consideration of where and when to recharge. 2) Flexible platoon formation: aligning departure times and route choices for multiple trucks to form platoons on shared road segments. 3) Battery balancing: allowing en‑route platoon position swaps so that energy consumption is distributed wisely across vehicles, preventing any single truck from depleting its battery too quickly. 4) Labor cost integration: accounting for driver waiting time during charging stops, which directly impacts total driver labor expenses. To address these challenges, we construct an MILP that jointly determines routing paths, charging station selections, labor hours, and platoon position swaps. Exact solution methods are likely to be intractable for realistic networks, especially when the problem size scales, so we develop an Adaptive Large Neighborhood Search (ALNS) algorithm, enhanced by the savings bounds method, infeasible‑pair elimination as the preprocessing approach, and a candidate‑station filtering wrapper method to speed the solution convergence. Computational experiments on test instances with up to 150 trucks demonstrate that incorporating platooning can reduce total operational costs by up to 2.77\%. The proposed algorithm also shows strong computational efficiency, achieving up to 99.96\% savings in computation time compared to CPLEX and solving instances with 150 trucks in approximately 120 seconds, highlighting its potential for real-world applications.

\hfill\break%

\end{abstract}

\newpage
\section{Introduction}

Truck platooning has emerged as a promising logistics innovation for both conventional and emerging truck fleets. It is typically defined as the linking of two or more semitrailer trucks into a convoy that travels at relatively high speeds with short inter-vehicle headways. In this configuration, following vehicles experience reduced aerodynamic drag, which in turn improves overall energy efficiency and can substantially lower operating costs.

The advantages of truck platooning can be identified along several dimensions. First, it addresses the urgent need to reduce fuel consumption and greenhouse gas (GHG) emissions in the trucking sector. Studies such as \citep{platoonsaving, bishop2020} have shown that truck platooning can achieve fuel savings of up to 10\% for trailing vehicles and 4.5\% for leading vehicles, and that an average energy reduction of 8.48\% per truck is attainable when the maximum platoon size is seven. By substantially cutting fuel consumption, platooning contributes directly to climate change mitigation through lower emissions. Second, these fuel savings translate into significant economic benefits. The \cite{atri} estimates that fuel accounts for roughly 24\% of a trucking company's operating cost in 2024; therefore, a 10\% reduction in fuel consumption corresponds to approximately 2.4\% total cost savings. Third, truck platooning can improve road safety. With synchronized acceleration and braking supported by wireless communication, trucks in a platoon are less prone to accidents caused by sudden stops or human error. Moreover, by reducing the manual workload of vehicle control—particularly for following trucks—platooning has the potential to alleviate driver fatigue \citep{larsen2019hub,sun2021auction, sun2021decentralized, Janssen2015TruckPD}. As a result, fewer accidents, smoother traffic flow, and more reliable travel times can be expected, enhancing the overall efficiency of road transport. From a macroscopic perspective, these improvements support a more stable and resilient supply chain \citep{supplychain}. Furthermore, by partially automating longitudinal control, truck platooning may eventually enable reductions in labor costs and even the number of drivers required, which represent a substantial share of logistics operating expenses.

Given these benefits—including substantial energy savings, reduced emissions, enhanced safety, and potential labor efficiencies—truck platooning is widely recognized as a transformative technology for the trucking industry. Consequently, a large body of research has focused on optimizing technical aspects such as truck speeds, longitudinal spacing, and route coordination to maximize platoon formation opportunities and operational efficiency. Advanced algorithms and vehicle-to-vehicle (V2V) communication technologies have been explored to ensure precise coordination among trucks, with the overarching goal of minimizing aerodynamic drag and improving fuel efficiency. While these technical advancements are essential, fully realizing the potential of truck platooning requires integrating these capabilities into broader logistics and supply chain decision-making. In particular, it is crucial to understand how platooning interacts with real-world freight operations, including routing, scheduling, and resource allocation, under different logistics scenarios.

Within this broader context, energy consumption is one of the most critical dimensions where truck platooning can offer additional value. For conventional diesel fleets, fuel prices, although volatile, are relatively well understood, and the marginal benefits of fuel savings are meaningful but somewhat predictable. In contrast, for electric trucks, the sensitivity to energy efficiency and charging costs is much higher, and this sensitivity is becoming more pronounced as electric trucks gain market share. Thanks to advances in battery technology, many electric trucks are now advertised with per-charge ranges of approximately 150–400 miles \citep{nacfe2024}, which align well with typical distances between regional distribution centers and fulfillment hubs. This makes them particularly suitable for mid-haul operations, where overnight charging at warehouses can minimize downtime and mitigate range anxiety. However, real-world factors such as battery degradation and adverse weather conditions can significantly reduce effective driving range \citep{nacfe2024}, introducing greater operational uncertainty than is typically encountered with diesel vehicles.

Truck platooning directly addresses this energy challenge by lowering energy consumption and extending the effective range of electric trucks. For electric fleets, this yields several critical advantages. First, charging prices exhibit substantial temporal and spatial variability across stations. By reducing energy consumption, platooning gives trucks more flexibility to bypass expensive charging locations and instead recharge at cheaper stations, thereby lowering overall operating costs. Second, an extended effective range reduces both the frequency and duration of charging stops. Fewer stops improve delivery reliability and reduce the total cost of electricity and associated labor time at charging facilities. Additionally, platooning enables the redistribution of State of Charge (SOC) among trucks by alternating the leading position. Since the leading truck experiences higher energy consumption due to increased aerodynamic drag, rotating this position among vehicles balances energy usage across the platoon. This strategy not only extends the collective range of the platoon but also prevents any single truck from depleting its battery disproportionately. As the logistics industry accelerates its transition toward electrification, the combination of platooning and electric trucks thus offers a particularly compelling pathway: reduced energy consumption, optimized charging strategies, and balanced in-platoon energy usage jointly support greenhouse gas reduction targets while enhancing cost-efficiency and service reliability. Together, these capabilities position electric truck platooning as a keystone for sustainable and scalable freight operations.

Despite this potential, most existing studies on the operation planning of electric truck platooning have focused on simplified highway-corridor settings. Such models only partially capture the realities of regional trucking logistics, where electric trucks operate over complex road networks and must strategically plan intermediate charging stops to manage range limitations. In practice, delivery trucks rarely travel along a single highway; instead, they traverse multiple interconnected routes on a road network, facing diverse options for path and charging-station selection.

Motivated by these gaps, this paper studies electric truck platooning in a more realistic network setting and makes the following contributions. First, to the best of our knowledge, we provide the first optimization framework for planning electric truck operations on a general road network rather than along a single corridor. In this setting, candidate charging locations for each truck are not confined to one path: trucks can choose among multiple routes between their origin–destination pairs, inducing a combinatorially large space of routing and charging plans. Second, we incorporate both charging costs and drivers’ labor costs into the operational cost function to better reflect the decision process of a central coordinator. In particular, detours taken to form platoons may reduce energy cost but increase labor hours, introducing an explicit trade-off between platoon savings and labor expenses. Third, we allow trucks in a platoon to swap positions en route, so that they can alternately take the lead position and balance their energy consumption, thereby extending the effective travel range of the platoon as a whole.

The remainder of this paper is organized as follows. Relevant studies on electric truck platooning and related planning models are reviewed in Section~\ref{sec:literature}. The problem is mathematically formulated in Section~\ref{sec:formulation}, together with an illustrative example. A novel iterative solution algorithm is presented in Section~\ref{sec:alg}. In Section~\ref{sec:exp}, numerical experiments are conducted to demonstrate the performance of the proposed modeling framework and solution algorithm, as well as to analyze the impact of several key parameters on system costs. Finally, Section~\ref{sec:con} concludes the paper and outlines directions for future research.

\section{Literature} \label{sec:literature}

As transportation electrification has grown significantly in recent years, electric trucks have also gained public attention as a strategy to build greener and environmentally friendly logistics. Studies such as \cite{OSIECZKO2021113177} and \cite{IMRE2021102794} have highlighted electric trucks' potential to reduce greenhouse gas emissions in freight operations. In real-world settings, \cite{boriboon} analyzed real-world activity patterns of heavy-duty battery electric trucks in regional distribution fleets, finding that electric trucks not only substantially lower emissions but may also meet the range requirements typically necessary for daily freight operations. However, the limited range of electric trucks remains a key challenge, prompting many early studies such as \cite{BASSO2019141, ANDERLUH2021940} to focus on last-mile logistics, often overlooking recharging needs. More recent efforts, such as \cite{MACRINA2019183}, have started to incorporate charging behavior to address longer delivery routes more accurately.

Despite the natural synergy between platooning and electrification, where energy savings through reduced aerodynamic drag can alleviate range anxiety, studies that simultaneously investigate charging decisions and platooning coordination remain scarce. A few recent works have attempted to bridge this gap. For example, \cite{eplatoonmanage} proposed an energy management framework for electric truck fleets, emphasizing platooning, cargo allocation, and charging decisions across three operational scenarios: pre-dispatch cargo allocation, sequence switching during highway travel, and platooning coordination at charging stations. The study also presented closed-form analytical solutions and LP-based heuristics to maximize travel range and minimize fleet charging time, with an innovative focus on platoon sequence switching. However, the author assumes a simplified highway-only problem instead of a road network, overlooking the routing and scheduling under a more realistic setting.

Similarly, \cite{ALAM2023104009} formulated a co-optimization problem integrating electric vehicle charging and platooning using a branch-and-cut algorithm on a 595-mile Florida highway corridor. Their model incorporated detailed constraints, including battery dynamics, charging station availability, and platoon coordination, providing strong computational insights into the routing and energy planning of mid to long-haul electric trucks. In another recent work, \cite{eplatoon} developed a metaheuristic search framework to jointly optimize routing, charging, and platooning strategies for electric trucks. Their model considers dynamic platoon formation opportunities at key intersections, charging cost minimization, and vehicle assignment. While these studies are comprehensive in scope, they do not consider flexible driver dispatch or labor costs, which are critical for cost optimization in real-world freight networks. In addition, they restrict the problem context to a highway corridor instead of road networks, failing to exploit the platoon's benefit under more realistic scenarios. 

In addition, a recent study by \cite{sun2025} is highly relevant and proposed a highway-network-based routing model for electric truck platooning that simultaneously considers driver assignments, charging decisions, and platoon coordination. They formulated an MILP to address these interdependent components and provided a Lagrangian Relaxation approach to solve the problem efficiently. While their work is among the most detailed in jointly modeling charging and platooning, it still assumes static platoon leadership and does not consider energy rebalancing through position rotation within the platoon. Another relevant study by \cite{FuChow2024TRB} introduced a modular, platoon-based vehicle-to-vehicle electric charging framework that explicitly models energy transfers inside a platoon and on a road network. However, the focus of this paper was to minimize energy consumption and travel time instead of minimizing the operation cost of a logistics provider. They also did not take the energy consumption savings of following trucks as well as variable charging costs into consideration.

\section{Problem statement} \label{sec:formulation}

In this section, we examine the impact and benefits of truck platooning for mid-haul to long-haul logistics powered by electric truck fleets, especially on how the energy savings brought by truck platooning may alter the charging consideration of electric trucks. Similar to private electric vehicles, electric trucks encounter challenges such as range anxiety and delays due to frequent recharging. Truck platooning offers a potential solution by reducing energy consumption, alleviating range anxiety, and improving overall trucking efficiency. To better reflect real-world scenarios, the problem is modeled within the context of a road network, incorporating en-route charging stations. This study will focus on how a carrier company can dispatch its electric truck fleet to meet delivery demands while leveraging platooning to minimize total charging and labor costs. Additionally, we will explore the potential effects of electric truck platooning on decision-making processes, including the selection of charging stations, as well as its impact on routing and scheduling strategies for mid- to long-haul logistics operations.

\subsection{Illustrative example}\label{sec:descrip}
\begin{figure}[htbp]
	\centering
         \subfigure[Example network]{\label{fig:evsamplenet}\includegraphics[scale=0.4]{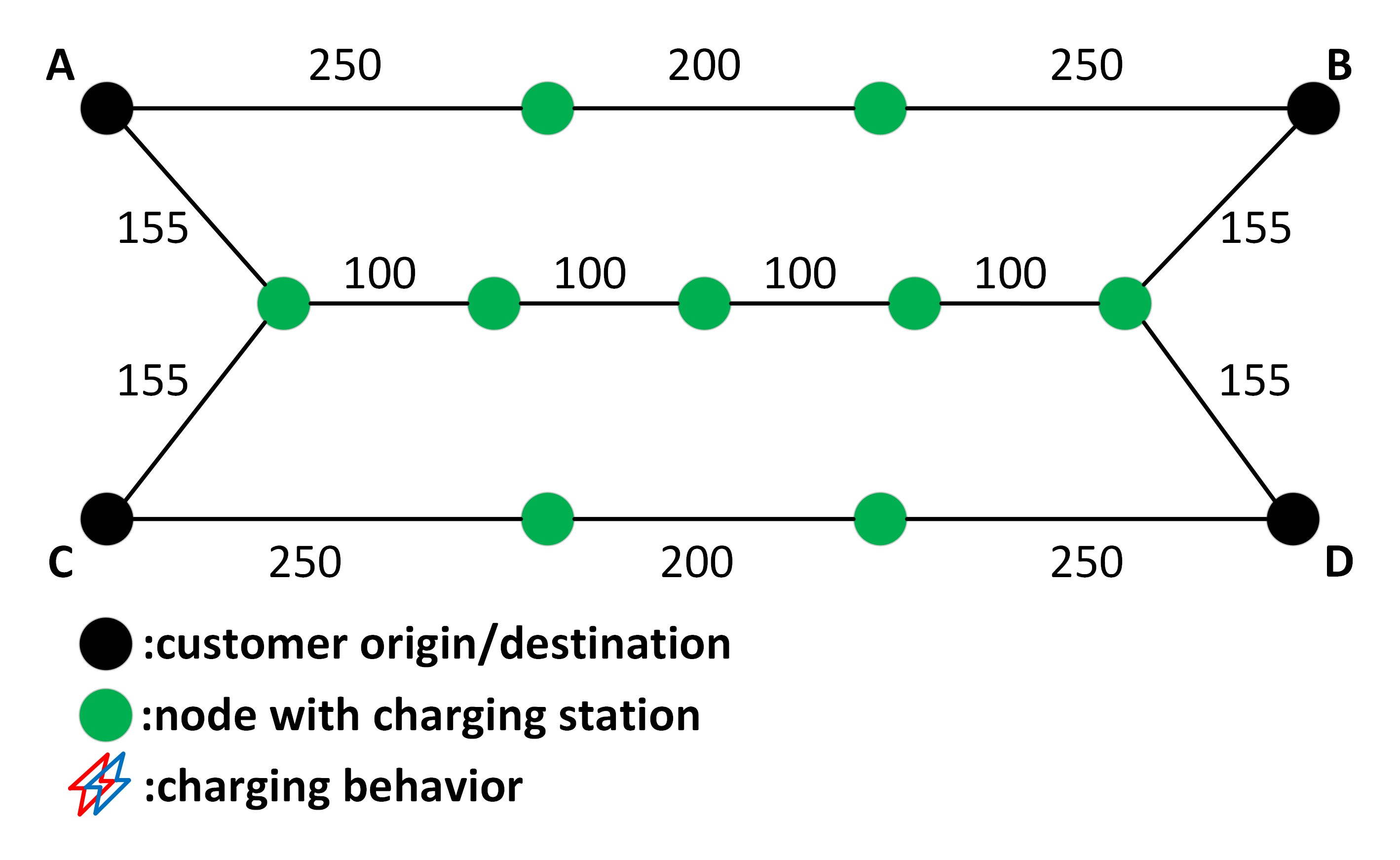}}\quad
	\subfigure[Non-platooning case]{\label{fig:evsamplenop}\includegraphics[scale=0.4]{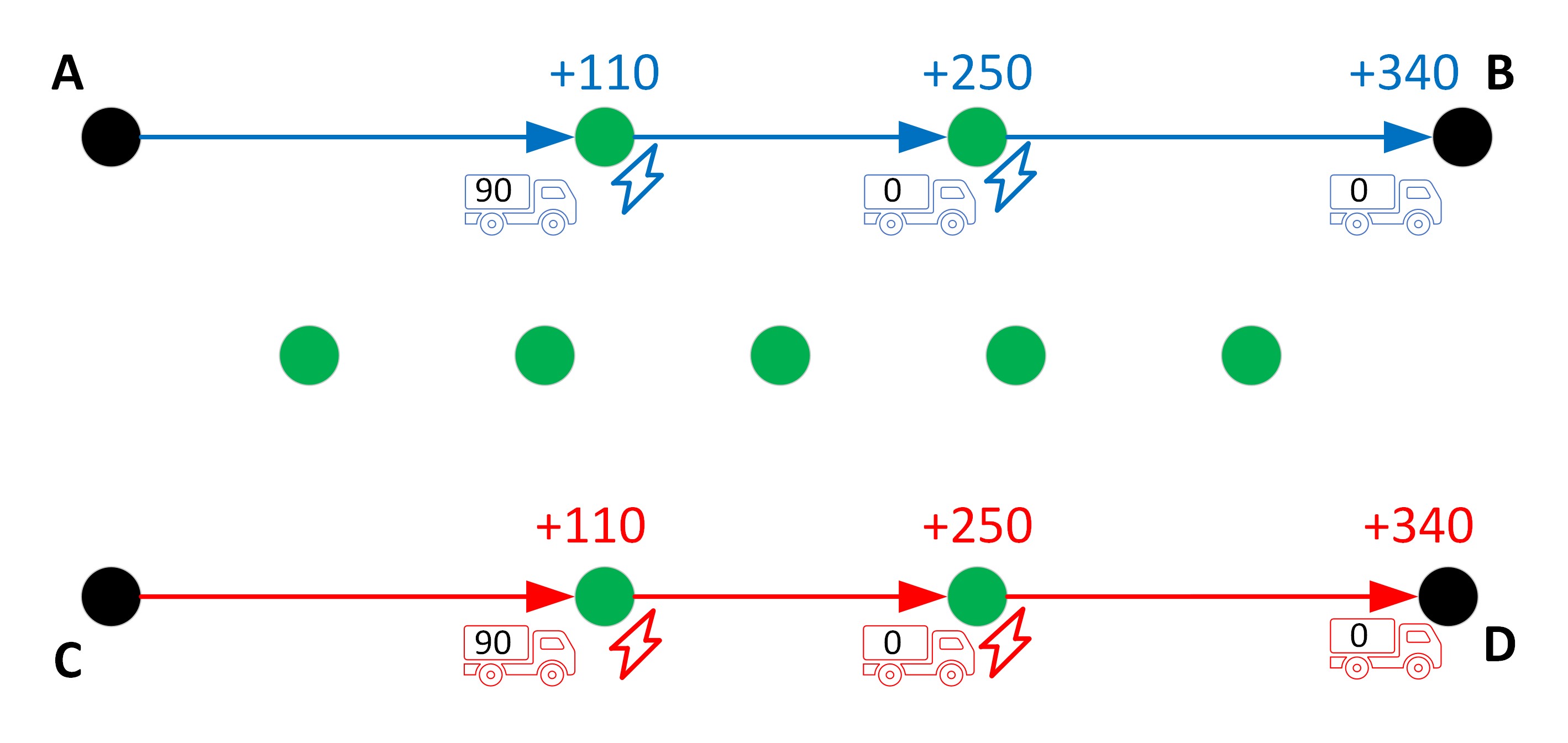}}\quad
    \subfigure[Platooning case]{\label{fig:evsampleplatoon}\includegraphics[scale=0.4]{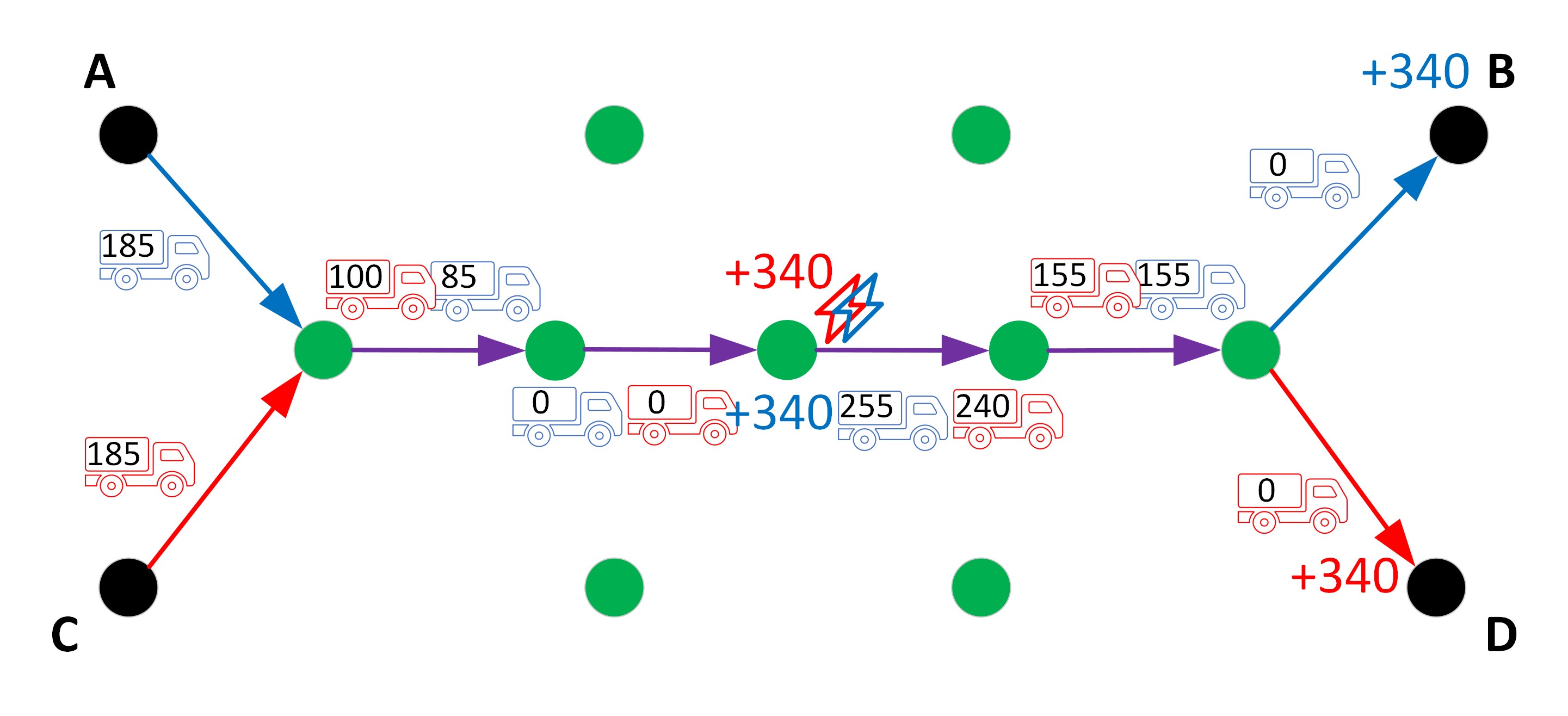}}\quad
    \caption{Illustrative example}
    \label{fig::platoonelectric}
\end{figure}

As mentioned before, electric trucks face challenges such as range anxiety and frequent recharging. To illustrate the potential benefits of electric truck platooning in mitigating range anxiety and optimizing charging strategies, we present a simple example based on the network shown in Figure \ref{fig:evsamplenet}. The network consists of two origin nodes ($A$, $C$), two destination nodes ($B$, $D$), and several intermediate nodes connected by highways with distances labeled in kilometers, and all green nodes represent the charging stations. To serve the customer demand, one electric truck is dispatched from $A$ to $B$, and another truck is dispatched from $C$ to $D$.

The key assumptions for this example are as follows:

\begin{itemize}
    \item The maximum driving range of a fully charged electric truck is 340 km, and each truck starts its trip with a full battery and should be recharged fully after finishing its trip.
    \item The speed is identical en route and yields a constant energy consumption rate, which is assumed to be 0.2 kWh/km.
    \item Charging price and power are assumed to be identical across different charging stations, and we assume the charging price is \$0.5/kWh.
    \item Labor cost is omitted for simplicity.
    \item Platoons can be formed at any intermediate nodes.
    \item The platooning energy-saving ratio is set at 15\%  to account for aerodynamic drag reduction and regenerative braking benefits specific to electric vehicles \citep{peloton, bishop2020}.
\end{itemize}

Two scenarios are considered: (1) trucks operate independently without platooning, and (2) trucks cooperate to form a platoon along shared segments.

In the non-platooning case, each truck travels its shortest path and needs to make two stops to charge in order to reach destination. The minimum possible mileage charged at the intermediate charging stations is shown in Figure \ref{fig:evsamplenop}, with lightning signs of different colors on the side to represent the charged amount for different trucks. The numbers represent the recharged driving range with units of km instead of kWh to better present how much driving range has been restored for each truck. Furthermore, the numbers inside each truck specify the remaining driving range after reaching the end of the link. For example, the truck departs from A with full capacity will have a remaining driving range of 90 upon arrival at the first charging station, as shown in Figure \ref{fig:evsamplenop}. Accordingly, in the non-platooning case, the total charging cost for trips from $A$ to $B$ and $C$ to $D$ (via intermediate nodes) is computed as: $(110 + 250+340) \times 0.2 \times 2 \times 0.5 = \$140$.

In the platooning case, the trucks synchronize their schedules to travel together on the common path in between. Specifically, due to platoon savings, the following truck now only consumes 85\% of energy on the common path. It is shown in Figure \ref{fig:evsampleplatoon} that the red truck first acts as a following truck after merging with the blue truck as a platoon, so it saves 15km mileage in the first common segment, and then the blue truck takes the following position and is managed to use the rest of 85km mileage to travel the 100km segment. Therefore, both trucks have enough energy to arrive at the intermediate charging station before running out of battery. They then recharge fully for the 340km driving range. After the recharge, the red truck first acts as the leader for the first road segment, and then the blue truck acts as the leader for the latter road segment. They can then arrive at their corresponding destinations without additional stops. Under the platooning case, we can witness that trucks only need to stop and charge once to complete their trips. In this case, the total energy cost is:$ (340+340) \times 0.2 \times 2 \times 0.5 = \$136$. Therefore, platooning provides a total $\frac{140 - 136}{140} = 2.86\%$ cost savings by energy reduction. 

This simple example provides an initial insight on the benefits that platooning can bring for the mid-haul electric truck logistics. Based on the example shown, with platooning benefit, both trucks are able to travel a longer distance, so they can take an alternative route to platoon together instead of traveling on their shortest path. Despite each traveling a slight longer distance to fulfill its trip, it saves energy consumption, which in turn saves the charging cost, and it also make fewer stops en route which may benefit a lot more when we take labor hours into consideration. 

As we further analyzed from Figure \ref{fig:evsampleplatoon}, we can see that platoon formation changes at every intermediate node in order to balance the energy consumption between the two trucks. Otherwise, if the second and the fourth intermediate nodes are removed from the network, leaving a 200-distance link connecting the adjacent nodes, then two trucks will not be able to have platoon reformation after traveling 100km, and thus one of the trucks will definitely consume more energy than the other. Therefore, it ends up that the leader needs 200km of driving range of energy, and the follower needs only 170km, which means that the follower does not have enough energy to reach the middle charging station without stop for charging once more. Therefore, in our problem statement, we consider the leader swapping feature during traveling, which allows different trucks in the same platoon to travel as a leader for different times, providing possibility for the trucks in the same platoon to balance their energy consumption. In this specific example, if the leader swapping is enabled, even if the two nodes are removed, both trucks can act as a leader for 50\% of the time when traveling as a platoon, and each will consume the energy of 185km drive range, so they are still be able to arrive at the middle charging station without further stops.

This simple example demonstrates how platooning can alter not only energy consumption but also route choices and charging strategies. It further elaborates the motivation to consider leader swapping in the problem context and our following mathematical formulation because of its utility in balancing energy use and enhancing operational efficiency.

\subsection{Mathematical formulation}\label{sec:form}

Generally, in this problem, we want to investigate the impacts and benefits of truck platooning under the electric truck delivery context. Therefore, the problem will be formulated as a pickup and delivery problem with truck platooning enhancement over a road network under mid-haul logistics background. Accordingly, there are several key assumptions need to be addressed: 1) Trucks have pre-determined origin and destination in the pickup-delivery problem setting. 2) We only include latest arrival times for each truck's time window, without earliest departure time. 3) Trucks' battery capacity, traveling speed, and energy consumption are assumed to be identical for simplicity without the loss of generality.  4) Charging speed across charging stations are assumed to be identical, but charging prices may vary. 5) Platoon saving ratio for energy consumption is the same no matter the platoon size. 6) We allow leader swapping in the same platoon while traveling to balance energy consumption. 7) Platoon formation and disintegration times are neglected, but there will be associated cost to reflect labor effort. 8) Every truck departs from its origin with battery capacity as $SOC^u$, and it is required to recharge to $SOC^u$, the maximum allowable state of charge after arriving at destination. With all the above setting and assumptions, we stand as a point of view of a central coordinator to optimize the routing, scheduling, and charging plans of trucks and platoons to minimize the total operation cost, which consists of charging cost, labor cost, and platoon restructuring/formation costs. 

Before presenting the overall model, the notations, including parameters and decision variables, frequently used in the mathematical model are provided in Table \ref{list_of_notations_2}.

\begin{table}[htbp]
\caption{List of Notations}
\label{list_of_notations_2}
\centering
\resizebox{\textwidth}{!}{
\begin{tabular}{cl}
\toprule
Notation                               & \multicolumn{1}{c}{Definition}\\ 
\midrule
\multicolumn{2}{c}{\bf Sets}\\ \hline
$P$ & Set of trucks' pickup positions\\
$D$ & Set of trucks' destinations\\
$C$ & Charging stations\\
$N$ & Set of all nodes, $N = P\cup D \cup C\cup O$\\
$A$ & Set of all arcs\\
$G$ & The entire graph of the road network, $G = (N,A)$\\
$K$ & Set of demand electric trucks\\
\hline
\multicolumn{2}{c}{\bf Parameters}\\ \hline
$L$ & Maximum platoon size\\
$t_{i,j}$ & Traveling cost on arc $(i,j)\in A$\\
$\tau^{LA}_{i}$ & Latest arrival time for truck $k\in K$, $i = d(k)$\\
$p(k)$ & Pickup location of truck $k, k\in K$\\
$d(k)$ & Destination location of truck $k, k\in K$\\
$\alpha_1$ & Unit in-vehicle labor cost with respect to leading truck time\\
$\alpha_2$ & Unit in-vehicle labor cost with respect to following truck time\\
$\alpha_3$ & Unit idling labor cost with respect to rest/dwell time at charging stations\\
$\alpha_4$ & Additional cost per platoon restructuring\\
$\beta$ & The platoon saving factor\\
$\sigma$ & The energy consumption rate with respect to time\\
$\eta$ & Charging speed of charging stations\\
$\gamma_i$ &The charging price per unit of electricity at charging station $i,i\in C$\\
$Q$ & The amount of electricity of a full-capacity battery of the electric truck\\
$SOC^u$ & The maximum allowable SOC \\
$SOC^l$ & The lowest allowable SOC before charging\\
\hline
\multicolumn{2}{c}{\bf Variables}\\ \hline
$x_{i, j, k}$ & Binary variables to indicate whether truck $k$ travels on arc $(i,j)$\\
$y_{i,k}$ & The SOC level of truck $k$ when it reaches node $i$\\
$l_{i, j, k}$ & Binary variable: 1 if truck $k$ is the leading vehicle of a platoon on arc $(i,j)$; and 0, \\
& otherwise. \\
$h_{i,j,k}$ & The percentage that truck $k$ acts as a leading truck, the range is between 0 and 1.\\
$f_{i, j, k_1, k_2}$ & Binary variable: 1 if truck $k_1$ is the leading vehicle and $k_2$ is in the same platoon on \\
& arc $(i,j)$; and 0, otherwise.\\
$v_{i, k}$ & Linear variable: the time truck $k$ charges at node $i\in N$, 0 if not charged.\\
$s_{i,k}$ & Arrival time of truck $k$ at node $i$\\
$w_{i, k}$ & The dwell time for truck $k$ at node $i$\\
$\tilde{h}_{i,j,k_1,k_2}$ & Dummy variable with a range from 0 to 1 to represent the product of $f_{i,j,k_1,k_2}$ and $h_{i,j,k_2}$.\\
$e_{i,j,k}$ & A binary dummy variable to indicate whether truck $k$ acts as a leading  truck on arc $(i,j)$\\
\bottomrule
\end{tabular}}
\end{table}

\subsubsection{Objective function}

\begin{flalign}
\min z = & \eta \sum_{i\in N}\sum_{k\in K}\gamma_i  v_{i,k} + \sum_{(i,j)\in A} t_{i, j}\bigg[\alpha_1\sum_{k\in K}l_{i, j, k}+\alpha_2\sum_{k_1\in K}\sum_{k\in K:k \neq k_1} f_{i, j, k_1, k}\bigg] &\notag\\
&+\alpha_3\sum_{i\in N\setminus\{p(k),d(k)\}}\sum_{k\in K}w_{i,k}  +\alpha_4 \sum_{(i,j)\in A} \sum_{k\in K} (e_{i,j,k} - l_{i,j,k}) &\label{obj2}
\end{flalign}
where $v_{i,k}$ is a variable to specify the amount of electricity charged for truck $k$ at node $i$ with charging station; binary variable $l_{i,j,k}$ tells if truck $k$ acts as a leading vehicle of a platoon on link $(i,j)$, and $f_{i,j,k_1,k}$ indicates if truck $k$ is following $k_1$ on link $(i,j)$; $W_{i,k}$ portraits the dwell time of truck $k$ at node $i$; $p(k)$ and $d(k)$ specify the origin and destination of each demand truck $k$; and $e_{i,j,k}$ is a binary variable to indicate if this truck $k$ once acts as a leader while traveling on link $(i,j)$. Obviously, any truck can only act as either a leading truck or a following truck at any time, and there is only one leading truck for any platoon. 

The objective minimizes the total operation cost of the central coordinator, which can often be a carrier company, including the charging cost and labor cost. Specifically, the first term $\eta \sum_{i\in N}\sum_{k\in K}\gamma_i  v_{i,k}$ is the charging cost, and it is obtained by multiplying the total charging amount $\eta v_{i,k}$ by the unit charging price $\gamma_i$ at each station, varying from station to station. The labor cost is composed of three parts: the in-vehicle delivery cost, the idling cost, and the platoon restructuring cost. The second term refers to the in-vehicle delivery labor cost, while the third one refers to the idling cost. Drivers of the leading vehicles should be paid more, while the drivers of the following vehicles are paid less for less workload. In addition, drivers who have rested during dwell time will only take the lowest wage. Therefore, we assign $\alpha_1$, $\alpha_2$, and $\alpha_3$ with different values to depict the wage level differences. The last term represents the additional cost for any platoon restructuring during traveling. Whenever an original follower in the platoon changes to a leader during traveling, we identify it as a leader swapping behavior and then additional cost will be incurred. $\alpha_4$ is the associated weight coefficient used to depict the additional effort paid for the platoon restructuring.  

\subsubsection{Truck flow constraints}

\begin{flalign}
&\sum_{i\in N,(i,j)\in A} x_{i, j, k}=\sum_{h\in N, (j,h)\in A} x_{j, h, k},
&\forall j\in N\setminus \{p(k),d(k)\}, k\in K, & \label{et2}\\
&\sum_{j\in N,(i,j)\in A} x_{i, j, k}=1,
&\forall i = p(k), k\in K,&\label{et3}\\
&\sum_{i \in N, (i,j)\in A} x_{i, j, k} =  1
&j = d(k), \forall k\in K, &\label{et4}
\end{flalign}
where $x_{i,j,k}$ tells whether truck $k$ travels on link $(i,j)$.

Constraints (\ref{et2})-(\ref{et4}) construct the routes of trucks and ensure that each truck must depart from its origin and arrive at its destination, and flow conservation at each intermediate node visited. Constraints (\ref{et2}) guarantee that the inflow of trucks should be equal to the outflow of trucks at any node, except for the origin and the destination. Constraints (\ref{et3}) and (\ref{et4}) ensure that each truck must depart from its origin and arrive at its destination.

\subsubsection{SOC constraints}

\begin{flalign}
&y_{i,k} \geq SOC^l
&\forall i\in N, k\in K, &\label{et5}\\
&y_{i,k}\cdot Q + \eta \cdot v_{i,k} \leq SOC^u\cdot Q
&\forall i\in N, k\in K, &\label{et6}\\
&y_{i,k} = SOC^u
&i = p(k), \forall k\in K,  &\label{et7}\\
&y_{i,k}\cdot Q + \eta \cdot v_{i,k} = SOC^u\cdot Q
&i = d(k), \forall k\in K, &\label{et8}\\
& v_{i,k} = 0
&i\in \{N \setminus C\}\cup \{p(k)\}, \forall k\in K,  &\label{et9}
\end{flalign}
\vspace{-10mm}
\begin{flalign}
&y_{i,k}\cdot Q + \eta \cdot v_{i,k} - \sigma t_{i,j}\bigg[h_{i,j,k} + (1-\beta)(1-h_{i,j,k})\bigg] + M(1-x_{i,j,k}) \geq y_{j,k}\cdot Q 
&  &\notag \\ 
&\forall (i,j)\in A, k\in K,  &\label{et10}    
\end{flalign}

\noindent where $y_{i,k}$ records the remaining battery state of charge (SOC) of each electric truck $k$ reaches node $i$, and we assume the maximum capacity, $Q$, of each truck's battery is identical. $SOC^l$ and $SOC^u$ are chosen to be the least and the most allowable SOC. $\beta$ is the platoon saving ratio on electricity consumption for the following trucks. $\sigma$ represents the energy consumption rate with respect to time, and $h_{i,j,k}$ indicates how much percentage of time a truck acts as the leader in a platoon while traveling. $M$ is a sufficiently big number.

Constraints (\ref{et5}) and (\ref{et6}) illustrate that the battery amount of each truck should not fall below the least allowable SOC, while also not exceeding the most allowable SOC at any time, no matter in charging or not. Constraints (\ref{et7}) tell that every truck starts its trip with the maximum allowable battery amount. Constraints (\ref{et8}) guarantee that every truck should be charged to the original SOC when it departs, and we assume that every truck is always allowed to charge at its destination, but a truck cannot charge at its origin, as realized by constraints (\ref{et9}). Last but not least, constraints (\ref{et10}) track the SOC changes of each truck. If truck $k$ travels from node $i$ to $j$, the SOC when reaching node $j$ can be obtained by adding the charging amount to the SOC at node $i$, and then subtracted by the consumption amount while traveling on link $(i,j)$. 

\subsubsection{Time consistency constraints}

\begin{flalign}
&s_{i, k} + t_{i, j} + w_{i, k} - M(1 - x_{i, j, k}) \leq s_{j, k}
&\forall (i, j)\in A,  k\in K, &\label{et11}\\
&s_{i, k} + t_{i, j} + w_{i, k} + M(1 - x_{i, j, k}) \geq s_{j, k}
&\forall (i, j)\in A,  k\in K, &\label{et12}\\
&s_{i, k} \leq \tau^{LA}_{i} 
&i=d(k), \forall k\in K,  &\label{et13}\\
&w_{i,k} \geq v_{i,k}
&\forall i\in N, k\in K, &\label{et14}
\end{flalign}
where $s_{i,k}$ records the arrival time of truck $k$ at node $i$, and $w_{i,k}$ is the dwell time of truck $k$ at node $i$. 

Constraints (\ref{et11})-(\ref{et14}) honor the time consistency of each truck along its route. Constraints (\ref{et11})-(\ref{et12}) track the time schedules of each truck during its trip. Constraints (\ref{et13}) require each truck demand to be delivered to its destination before its latest arrival time. Constraints (\ref{et14}) indicate that the dwell time of each truck should be no less than its charging time at any specific node.

\subsubsection{Platoon constraints}

\begin{flalign}
&x_{i,j,k} = l_{i,j,k} + \sum_{k_1\in K,k_1\neq k}f_{i,j,k_1,k} 
&\forall (i,j)\in A, k\in K, &\label{et15}\\
&x_{i,j,k} - h_{i,j,k} \geq 0
&\forall (i,j)\in A, k\in K, &\label{et16}\\
&h_{i,j,k_1}+\sum_{k_2\in K, k_2\neq k_1}(f_{i,j,k_1,k_2} \cdot h_{i,j,k_2}) \geq 1 - M(1-l_{i,j,k_1})
&\forall (i,j)\in A, k_1\in K, &\label{et17}\\
&h_{i,j,k_1}+\sum_{k_2\in K, k_2\neq k_1}(f_{i,j,k_1,k_2} \cdot h_{i,j,k_2}) \leq 1 + M(1-l_{i,j,k_1})
&\forall (i,j)\in A, k_1\in K, &\label{et18}\\
&\sum_{k_2\in  K:k_2 \neq k_1} f_{i, j, k_1, k_2} \leq (L - 1) \times l_{i, j, k_1}
&\forall (i, j)\in A, k_1\in K, &\label{et19}
\end{flalign}
\vspace{-8mm}
\begin{flalign}
&-M(1 - f_{i, j, k_1, k_2}) \leq s_{i, k_1} + w_{i, k_1} - s_{i, k_2} - w_{i, k_2}
&\forall (i,j)\in A, k_1,k_2\in K, k_1\neq k_2, &\label{et20}\\
&s_{i, k_1} + w_{i, k_1} - s_{i, k_2} - w_{i, k_2} \leq M(1 - f_{i, j, k_1, k_2})
&\forall (i,j)\in A,  k_1,k_2\in K, k_1\neq k_2, &\label{et21}
\end{flalign}
\vspace{-8mm}
\begin{flalign}
&e_{i,j,k} \geq h_{i,j,k}
&\forall (i,j)\in A, k\in K, &\label{et22}
\end{flalign}
where $L$ is the platoon size limit.

Constraints (\ref{et15})-(\ref{et22}) honor the platooning consistency among trucks on the roads. Specifically, constraints (\ref{et15}) demonstrate that each truck can only be a leading truck or a following truck at the same time, and if it is traveling alone, it will be a leading truck of a platoon with only one vehicle. Constraints (\ref{et16}) illustrate that the percentage of a truck that can act as a leader on a link can only be positive when it does traverse such a link. Constraints (\ref{et17}) and (\ref{et18}) tell us that while we allow trucks in the same platoon to be the leading truck alternatively, the time ratio that each truck acts as a leading truck should be summed up to 1. Constraints (\ref{et19}) limit the maximum platoon size. Constraints (\ref{et20}) and (\ref{et21}) guarantee that if two trucks belong to the same platoon on link $(i,j)$, they should travel the link at the same time. Constraints (\ref{et22}) state that if a truck $k$ acts as a leading vehicle for a portion of time on the link (i,j), then the indicator variable $e_{i,j,k}$ should be 1. 

\subsubsection{Integral constraints}
\begin{flalign}
&e_{i,j,k}\in \left\{0,1\right\}
&\forall (i,j)\in A,k\in K, &\label{et23}\\
&x_{i,j,k}\in \left\{0,1\right\}
&\forall (i,j)\in A,k\in K, &\label{et24}\\
&0\leq y_{i, k}\leq 1
&\forall i\in N, k\in K, &\label{et25}\\
&v_{i,k}\geq 0
&\forall i\in N, k\in K, &\label{et26}\\
&l_{i, j,k}\in \left\{0,1\right\}
&\forall (i,j)\in A,k\in K, &\label{et27}\\
&f_{i,j,k_1,k_2}\in \left\{0,1\right\}
&\forall (i,j)\in A,k_1,k_2\in K,k_1\ne k_2, &\label{et28}\\
&0\leq h_{i,j,k} \leq 1
&\forall (i,j)\in A,k\in K, &\label{et29}\\
&s_{i, k}\geq 0
&\forall i\in N, k\in K, &\label{et30}\\
&w_{i, k}\geq 0
&\forall i\in N, k\in K, &\label{et31}\\
&0 \leq \Bar{h}_{i,j,k_1,k_2} \leq 1
&\forall (i,j)\in A,k_1,k_2\in K,k_1\ne k_2, &\label{et32}
\end{flalign}
where $\Bar{h}_{i,j,k_1,k_2}$ is a dummy variable that represents the product of $f_{i,j,k_1,k_2}$ and $h_{i,j,k_2}, \forall k_1, k_2\in K, k_1\neq k_2$. This variable is used to linearize the constraints (\ref{et17}) and (\ref{et18}).

Constraints (\ref{et23}) - (\ref{et32}) specify the ranges of all decision variables.

To linearize constraint (\ref{et17}) and (\ref{et18}), we can have:
\begin{flalign}
&h_{i,j,k_1}+\sum_{k_2\in K, k_2\neq k_1}\Bar{h}_{i,j,k_1,k_2} \geq 1 - M(1-l_{i,j,k_1})
&\forall (i,j)\in A, k_1\in K, &\label{et33}\\
&h_{i,j,k_1}+\sum_{k_2\in K, k_2\neq k_1}\Bar{h}_{i,j,k_1,k_2} \leq 1 + M(1-l_{i,j,k_1})
&\forall (i,j)\in A, k_1\in K, &\label{et34}\\
&\Bar{h}_{i,j,k_1,k_2} \leq h_{i,j,k_2}
&\forall (i,j)\in A, k_1, k_2\in K, k_1\neq k_2, &\label{et35}\\
&\Bar{h}_{i,j,k_1,k_2} \leq f_{i,j,k_1,k_2}
&\forall (i,j)\in A, k_1, k_2\in K, k_1\neq k_2, &\label{et36}\\
&\Bar{h}_{i,j,k_1,k_2} \geq h_{i,j,k_1} - (1 - f_{i,j,k_1,k_2})
&\forall (i,j)\in A, k_1, k_2\in K, k_1\neq k_2, &\label{et37}
\end{flalign}

Therefore, the overall model for the road-network electric truck platooning problem with alternative leading vehicle and charging consideration is summarized as follows:

\begin{flalign}
&\textbf{Minimize}& &\notag \\
&Obj = \eta \sum_{i\in N}\sum_{k\in K}\gamma_i  v_{i,k} + \sum_{(i,j)\in A} t_{i, j}\bigg[\alpha_1\sum_{k\in K}l_{i, j, k}+\alpha_2\sum_{k_1\in K}\sum_{k_2\in K:k_2 \neq k_1} f_{i, j, k_1, k_2}\bigg]  &\notag\\
&+\alpha_3\sum_{i\in N\setminus\{p(k),d(k)\}}\sum_{k\in K}w_{i,k} +\alpha_4 \sum_{(i,j)\in A} \sum_{k\in K} (e_{i,j,k} - l_{i,j,k}) &\notag \\
&\mbox{s.t.}\quad (\ref{et2})-(\ref{et16}) \quad \texttt{and} \quad (\ref{et19})-(\ref{et37})&\notag
\end{flalign}

\section{Solution Methodology}\label{sec:alg}

\begin{figure}[htbp]
	\centering
    \includegraphics[scale=0.4]{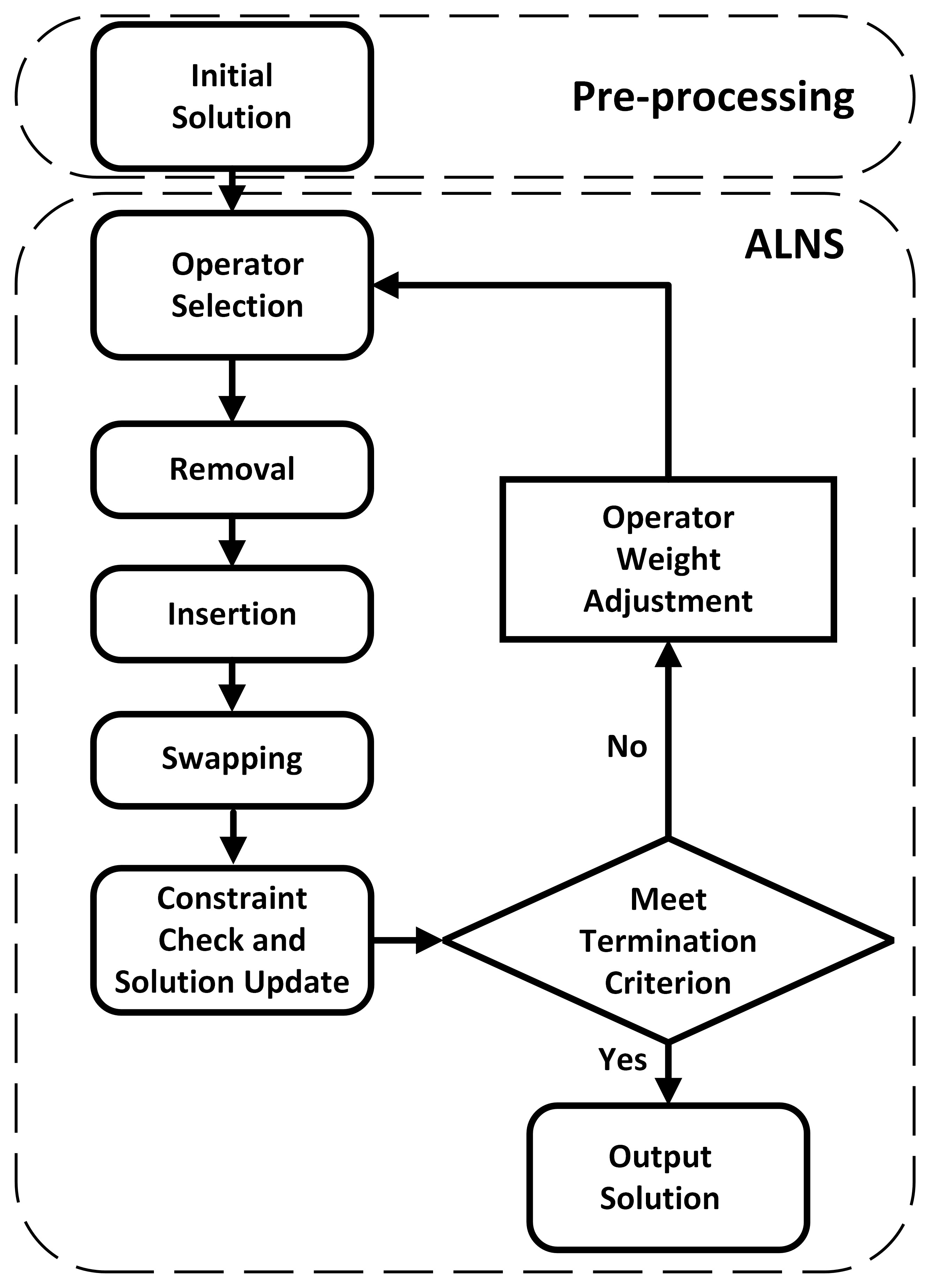}
	\caption{Flow chart of solution algorithm}
    \label{fig::evalg}
\end{figure}

The solution approach developed for the electric truck platooning problem comprises two main components: a pre-processing phase and an Adaptive Large Neighborhood Search (ALNS) algorithm, and the overall flow diagram is shown in Figure \ref{fig::evalg}. This design is motivated by the need to address the inherent complexity and large-scale nature of the problem, which involves the simultaneous optimization of routing, scheduling, charging, and platooning decisions.

ALNS was chosen as the core solution methodology due to its proven effectiveness in solving large-scale, combinatorially complex problems such as the Capacitated Vehicle Routing Problem (CVRP). Its flexibility in designing destroy and repair operators allows it to efficiently explore the vast solution space of electric truck platooning scenarios, which are further complicated by additional constraints like energy limitations, charging station availability, and dynamic platoon formations. 

The pre-processing phase plays a critical supporting role in this framework. It is designed to identify candidate charging stations for each truck, enabling the construction of a high-quality initial solution—an essential factor for the performance of ALNS. Moreover, the pre-processing step involves calculating both optimal and alternative paths for each truck to minimize combined charging and travel costs. These paths not only provide the foundation for initial solution construction but also supply key input to the regret-based insertion operators in ALNS. Specifically, the regret values used to prioritize delivery insertions are derived from the cost differences between optimal and sub-optimal paths identified during pre-processing.

This design effectively balances computational efficiency and solution quality. The pre-processing phase reduces the search space and guides ALNS toward promising regions of the solution landscape, while ALNS provides the adaptability and robustness needed to fine-tune solutions in the face of complex constraints and interactions among routing, charging, and platooning decisions.

\subsection{Pre-processing}
Pre-processing, shown as the first part of our algorithm in Figure \ref{fig::evalg}, includes several elements and provides an initial solution for the ALNS module. First, we examine each truck delivery and check if its optimal charging station places vary when it changes from traveling alone to traveling as a following truck. The intuition behind this is that when the truck is traveling alone, it cannot have any platoon benefits, so the truck will intend to take the path with minimum cost, charging cost and time cost; when the truck is always traveling as a following truck for its entire itinerary, it can obtain maximum platoon benefits. Accordingly, if the optimal charging stations are identical under such two circumstances, it means that such a truck will not likely alter its path with the existence of platoon potential, and thus these charging stations will more likely occur in the overall problem's optimal solution. 

However, we are not dealing with a simple shortest path problem, as charging costs vary across stations and en-route travel time still matters for the labor cost, so we developed an approach to determine which charging stations are favorable for each truck delivery, given the road network and charging prices for each station. To be noted, in our problem, we always assume that each truck needs to recharge to its original SOC, i.e., $SOC^u$, when arriving at the destination, and the charging price at its corresponding destination is always smaller than en-route charging stations. Since each truck may experience multiple charges along its route, we further assume that the maximum charging time except for the destination for each truck is 3 to mimic the real-world conditions.
Therefore, we need to consider a maximum of 4 possible scenarios for each truck: no charging, charging once, charging twice, and charging 3 times.

\subsubsection{Candidate charging station decision}
To decide a truck delivery's optimal charging plan, we will iterate from the no-charging scenario to the 3-charging scenario with 1 additional charge increment each time.

\begin{enumerate}
    \item Starting with the no-charging scenario, we use Dijkstra's algorithm ($O(n\log n)$ complexity) to obtain each truck's shortest travel path to minimize its total mileage. If the mileage of the path is smaller than the maximum SOC, the optimal charging plan is obtained, otherwise, go to the next step.
    \item When no-charging is infeasible, we consider charging once for each remaining truck delivery. To determine which charging station to adopt, it is equivalent to check the station $i$ that yields the minimum value of $\alpha_1(t(o,i) + t(i,d)) + \dfrac{(t(o,i) + t(i,d)-Q\cdot SOC^u)\cdot (\alpha_3+\gamma_i)}{\eta}$ with $t(o,i)$ and $t(i,d)$ representing the shortest travel time from origin to station $i$ and from station $i$ to destination. The running time of this check takes $O(n^3)$ complexity by employing Floyd-Warshall's algorithm. In addition, charging once will not be feasible whenever $t(o,i) + t(i,d)-2\cdot Q\cdot SOC^u > 0$ for all $i$.
    \item When charging once is infeasible, we consider charging twice for each remaining truck delivery, provided as follows:
    \begin{itemize}
  \item Preprocessing with Floyd Warshall Algorithm: Run Floyd Warshall's algorithm on the network to calculate the shortest paths and distances between all pairs of nodes. This step provides a foundation for efficiently querying the shortest distance between any two points, essential for subsequent optimization steps.
  
  \item Determining Energy Requirements: Given the final destination and the required minimum SOC to reach it, calculate the energy needed at the second charging station. This involves determining the shortest path from the second charging station to the destination and ensuring the truck has sufficient charge to complete this segment, factoring in the energy potentially remaining from earlier segments of the journey.
  
  \item Optimizing the Second Charging Station Selection: Search for all candidate second charging stations that matches the check: 1. $\min(t_{2nd, dest}) \leq Q\cdot SOC^u$ and 2. $\min(t_{2nd, origin}) \leq 2Q\cdot SOC^u$, where $\min(t_{2nd, dest})$ represents the shortest path travel time between the second charging station to destination, and $\min(t_{2nd, origin})$ as the shortest path travel time between the second charging station and the origin. The check is used to filter the stations that can fulfill the trip to the destination without additional stops, while also requiring only one additional stop to reach it from the origin. Then, we assign the required charging cost from $SOC^l\cdot Q$, along with idling cost and shortest path traveling cost, for each candidate second charging station so that the truck can fulfill the shortest path trip.
  
  \item Decomposing into a Single Charging Problem: With the calculated assigned cost at each candidate second charging station, the problem becomes analogous to finding a single optimal charging station for a journey from the origin to an intermediate point (the first charging station). This step essentially treats the segment from the origin to the first charging station as an independent problem, simplifying the optimization process. Then, for each candidate second charging station, we find its optimal first charging station and sum up the cost, and then we record the two smallest charging station pairs.
  
\end{itemize}
    \item When charging twice is still infeasible, we apply to charge 3 times, and it is similar to charging twice with an additional iteration.
\end{enumerate}

It is worth noting that this part can be further improved by formulating it as a discretized Markov Decision Process and solving it using Dynamic Programming to achieve results that are closer to optimal solution. The running time complexity can scale up to $O(n^3\log n)$.

\subsubsection{Initial Solution}
With the set of candidate charging stations for each truck determined, we partition the fleet into two groups based on the impact of platooning on their charging decisions: (1) a \textbf{no-difference group}, where platooning does not alter the candidate stations or charging plan; and (2) a \textbf{difference group}, where platooning leads to alternative charging station selections or timings.

The procedure begins by sorting all trucks according to their latest feasible departure times. These times are calculated as the latest possible arrival times at their respective destinations minus the travel time along their current sub-optimal candidate paths. Trucks are arranged in ascending order of these departure times to prioritize those with tighter scheduling constraints.

We then perform an iterative platoon formation process as follows:
\begin{enumerate}[label=\arabic*.]
    \item[\textbf{•}] Select the first truck from the no-difference group as the initial candidate for platoon formation.
    \item[\textbf{•}] Iteratively evaluate trucks from the difference group for possible insertion into the current platoon. Each truck is assessed based on two criteria: (\textit{i}) time feasibility—whether the truck can align its schedule with the platoon, and (\textit{ii}) energy or cost savings relative to traveling alone. Trucks satisfying both criteria are added to the platoon sequentially until no further feasible insertions are identified.
    \item[\textbf{•}] Proceed to the next truck in the no-difference group and repeat the process. First, check whether this truck can be integrated into any existing platoons formed in the previous steps. If not, initiate a new platoon following the same insertion logic.
\end{enumerate}

After all trucks in the no-difference group have been processed, any remaining trucks in the difference group are handled in a final iteration. Here, a single application of insertion operators is employed to integrate these trucks into existing platoons wherever feasible, or to form standalone routes if platooning is not advantageous.

\subsection{Adaptive Large Neighborhood Search}

With all the preprocessing explained above, we have the candidate charging station information for each truck as well as an initial solution that serves as a warm start input for the following ALNS method, which is showcased as the second module in the flow chart in Figure \ref{fig::evalg}. In addition, the optimal and alternative routing plans derived for each truck derived before will also be used to calculate the detour cost and regret value mentioned later in this section.

\subsubsection{Solution Encoding}
Each truck’s solution in the ALNS framework is represented in the following form:
$$
(1,1,0.1, (2,3),50\%)
$$
where each entry denotes:
\begin{itemize}
\item \textbf{Truck Index}: the unique identifier for the truck.
\item \textbf{Platoon Index}: the identifier of the platoon to which the truck belongs.
\item \textbf{Leading Ratio}: the proportion of distance or time for which the truck serves as the leading vehicle in the platoon.
\item \textbf{Links Traveled}: a list of network links traversed by the truck during its route.
\item \textbf{Charging Amount}: the amount of energy replenished at the latter node of the link, expressed as a percentage of full battery capacity. If the node is not a charging station, the amount must be 0. 
\end{itemize}

With each encoded truck segment representing the travel behavior of a truck on a specific link, the entire solution is constructed by the union of all encoded segments. In addition, based on the result of initial solution generation from the last subsection, we have the information about the route each truck takes and whether it platoons with other trucks. Given this information, we construct each truck's segments by decomposing its route to consecutive links, and assign the leading ratio to 0 or 1 based on platooning information. The platoon index is assigned iteratively until all truck segments have an associated platoon, and the charging amount is set as the minimum SOC to recharge to reach the next node.  

\subsubsection{Removal Operators}
\begin{enumerate}[label=\roman*.]
    \item \textbf{Worst truck delivery removal:}  
    This operator deletes $\kappa$ deliveries in non-decreasing order of their associated detour costs. Specifically, each truck's delivery is defined as the union of all the encoded truck segment of this truck, and the detour cost is defined as the difference between the current cost and the cost to travel the shortest path. The goal is to prioritize the removal of deliveries contributing least efficiently to the current routing plan. Specifically, the number of elements removed, denoted as $\kappa$, is drawn uniformly from a predefined range:  
\[
\kappa \sim U(\kappa_{\min}, \kappa_{\max})
\] 
where \(\kappa_{\min} = 2\) ensures at least minimal disruption because it guarantees to make some changes to the current solution that is composed of multiple encoded truck segments so that the solution has the potential to improve or jump out of local optima.  In addition, \(\kappa_{\max} = \max(2, \lfloor 0.1 \cdot K \rfloor\)) scales the upper bound of removal proportionally to the problem size \(K\), which is the number of trucks in our problem context. 

    \item \textbf{Random truck delivery removal:}  
    This operator deletes $k$ deliveries selected uniformly at random from the current solution, introducing diversification into the search process. 

    \item \textbf{Platoon non-exchange removal:}  
    The operator randomly selects a platoon for removal from the current solution. All deliveries assigned to this platoon are also removed. Only platoons that do not involve delivery exchanges with other platoons are considered. The probability of selecting a platoon for removal is inversely proportional to the number of deliveries it covers and the distance traveled as a platoon. 

    \item \textbf{Platoon removal:}  
    Similar to the previous operator, this one randomly selects a platoon for removal along with all its associated deliveries. However, platoons with delivery exchanges are also included in the candidate set. The selection probability is again inversely proportional to the number of deliveries covered, reflecting the platoon's relative efficiency.  

    \item \textbf{Exchanged delivery removal:}  
    This operator targets truck deliveries that traverse two or more different platoons, provided both platoons have a size greater than one. Deliveries with a higher number of exchanges are more likely to be selected for removal, as their removal can simplify the platoon coordination structure.  
\end{enumerate}

\subsubsection{Insertion Operators}

To reconstruct feasible solutions after removals, we design a set of insertion operators that incrementally reinsert deliveries into existing routes or platoons while considering energy consumption, charging, and time feasibility. The operators are described as follows:

\begin{enumerate}[label=\roman*.]
    \item \textbf{Greedy route insertion without a platoon:}  
    This operator inserts deliveries independently without considering platoon formation. Three selection criteria are used to prioritize deliveries:
    \begin{itemize}
        \item \textbf{Shortest path travel time:} Deliveries are sorted in descending order of their shortest path travel times from origin to destination.
        \item \textbf{Minimum charging and waiting cost:} Deliveries are prioritized based on the pre-processed minimum charging and waiting labor cost along their shortest path.
        \item \textbf{Regret value:} Deliveries are ranked by the absolute difference between the least and second least charging and waiting labor costs. Higher regret values indicate higher potential savings if inserted optimally.
    \end{itemize}

    \item \textbf{Greedy partial route insertion with a platoon:}  
    This operator attempts to form new platoons between pairs of truck deliveries. Given two deliveries, it identifies a pair of nodes $(a, b)$ such that both trucks can travel together on the link $(a, b)$, subject to a time feasibility check. Once a feasible shared link is determined, the operator will adjust those trucks' segments to have them platoon together on this specific link. The truck with more excess SOC at node b will be assigned as the leading truck. Furthermore,  deliveries are prioritized using the following criteria:
    \begin{itemize}
        \item \textbf{Shortest path travel time:} Deliveries with longer routes are prioritized to increase platooning potential.
        \item \textbf{Minimum charging cost with platooning:} Costs are calculated using a dynamic programming approach over a modified network, where links already traversed by trucks receive a cost discount if time feasibility allows for shared travel.
        \item \textbf{Regret value:} Similar to the previous operator, but calculated on the modified network, reflecting platooning and charging trade-offs.
    \end{itemize}

    \item \textbf{Partial route insertion into an existing platoon:}  
    This operator inserts a new delivery into an already established platoon by identifying overlapping links where the new truck can synchronize with the existing convoy. The procedure involves:
    \begin{itemize}
        \item Searching for candidate platoons where the truck’s route overlaps sufficiently with at least one link of the platoon.
        \item Evaluating time feasibility to ensure the truck can adjust its schedule for platoon formation without violating delivery constraints.
        \item Prioritizing insertions based on the expected additional savings in energy consumption and charging time due to platoon participation.
    \end{itemize}
    Deliveries are ranked using a hybrid score combining overlap length, potential energy savings from platooning, and incremental charging costs. This method encourages incremental expansion of existing platoons to improve overall efficiency.
\end{enumerate}

\subsubsection{Swapping Operators}

To further diversify the solution space and refine route configurations, two swapping operators are designed in this study:  

\begin{enumerate}[label=\roman*.]
    \item \textbf{Charging amount adjustment:}  
    This operator perturbs the charging amounts assigned to a truck at each charging station in its encoded solution segment.  The operator is applied subject to the SOC fulfillment check to ensure the truck maintains sufficient charge to traverse all planned route segments.  

    If a charging constraint violation occurs (e.g., insufficient SOC to reach the next planned station), the weight of this operator is dynamically increased to encourage more exploration of alternative charging amounts. For SOC adjustments:  
    \begin{itemize}
        \item If the SOC check indicates insufficient charge, the new charging amount is sampled uniformly from the range between \text{current amount} and $SOC^u - \text{current amount}$, allowing for an increase up to $SOC^u$.  
        \item Otherwise, if the SOC check indicates an overcharge, then the new charging amount is sampled uniformly from \([0, \text{current amount}]\), potentially reducing the charge to avoid unnecessary waiting time and charging cost at stations.  
    \end{itemize}

    \item \textbf{Random leading ratio adjustment:}  
    This operator modifies the proportion of time that a truck serves as the leading vehicle within a platoon. The adjustment is performed in increments or decrements of 0.1, ensuring that the total leading ratio across all trucks in the platoon sums to 1. This operator is only activated when the platoon size is two or more, as leading ratio adjustments are irrelevant for single-truck platoons. It helps balance energy consumption among trucks and supports equitable distribution of leading and following roles.
\end{enumerate}

\subsubsection{Constraint Checks}

To maintain feasibility after applying destroy and repair operators, a set of constraints is enforced:

\begin{enumerate}[label=\roman*.]
    \item \textbf{Leading ratio check:}  
    For each platoon, the sum of all trucks’ leading ratios must equal 1 to ensure a valid distribution of lead and follower roles.
    
    \item \textbf{Time feasibility check:}  
    Ensure that all trucks in a platoon can synchronize their schedules without creating infeasible time pairings. Specifically, each truck needs to meet its latest arrival constraint by taking charging time, dwell time, and trip time into consideration. 
    
    \item \textbf{Platoon size limit check:}  
    Verify that the number of trucks assigned to a platoon does not exceed the maximum allowable size specified by operational or regulatory constraints.
    
    \item \textbf{SOC fulfillment check:}  
    After swapping or adjusting charging assignments, confirm that the energy savings along each path segment (including between origin and destination or any two charging stations) are sufficient to meet the SOC requirements for the truck to complete its route.
\end{enumerate}

\subsubsection{Operator selection and acceptance criteria}

In the ALNS framework, destroy and repair operators are selected adaptively based on their historical performance. This is achieved by combining a softmax-based probability assignment (akin to a multinomial logit model) with a roulette wheel selection mechanism.  

At each iteration, the probability \(P_i\) of selecting operator \(i\) is calculated as:
\[
P_i = \frac{\exp(\lambda \cdot w_i)}{\sum_{j} \exp(\lambda \cdot w_j)}
\]
where:
\begin{itemize}
    \item \(w_i\) is the current weight of operator \(i\),
    \item \(\lambda = 10\) is a scaling parameter controlling the selection sensitivity (higher \(\lambda\) emphasizes differences in weights).
\end{itemize}

The roulette wheel mechanism is then applied to select an operator according to these probabilities. Specifically, the interval \([0, 1]\) is partitioned into segments proportional to each \(P_i\). A uniformly distributed random number \(u \in [0,1]\) is drawn, and the operator whose segment contains \(u\) is selected. This stochastic procedure ensures that higher-weight operators are favored while preserving a nonzero chance of selecting less frequently used operators, thereby promoting exploration.  

\textbf{Adaptive weight update mechanism:}

The weights \(w_i\) of each operator are updated adaptively every \(\nu = 50\) iterations. During each segment of \(\nu\) iterations, a performance score \(s_i\) is accumulated for operator \(i\). This score is incremented based on its contribution to solution improvements as follows:
\begin{itemize}
    \item If a new global best solution is found: \(s_i \leftarrow s_i + \delta_1\), where \(\delta_1 = 3\).
    \item If the solution improves relative to the current iteration: \(s_i \leftarrow s_i + \delta_2\), where \(\delta_2 = 2\).
    \item If a worse solution is accepted under the simulated annealing criterion: \(s_i \leftarrow s_i + \delta_3\), where \(\delta_3 = 1\).
    \item If the operator is not used in the segment: \(s_i\) remains unchanged.
\end{itemize}

At the end of each segment, the weights are updated using a smoothing scheme:
\[
w_i = (1 - r) \cdot w_i + r \cdot s_i
\]
where \(r = 0.2\) is the reaction factor determining the influence of recent performance.

\textbf{Simulated Annealing acceptance criterion:}

To allow the search to escape local optima, worse solutions may be accepted probabilistically based on a simulated annealing criterion. The acceptance probability for a worse solution is:
\[
P(\text{accept}) = \exp\left(-\frac{\Delta}{T}\right)
\]
where:
\begin{itemize}
    \item \(\Delta = f(s') - f(s)\) is the cost increase from the current solution \(s\) to the candidate solution \(s'\),
    \item \(T\) is the current temperature, initialized as \(T_0 = 100\) and updated by the cooling schedule:
    \[
    T_{k+1} = \alpha \cdot T_k
    \]
    where \(\alpha = 0.95\) is the cooling rate.
\end{itemize}

This combination of adaptive operator weighting, roulette wheel selection, and simulated annealing acceptance balances intensification and diversification, guiding the search effectively through complex solution spaces.

\section{Numerical Experiments}\label{sec:exp}

We assume that each electric truck has a constant battery capacity supporting a maximum driving range of 340\,km. According to a recent report \citep{evbattery}, the battery sizes of medium-duty trucks typically range from 92\,kWh to 139\,kWh. Without loss of generality, we adopt 135\,kWh as the representative battery size for this experiment. For simplicity and ease of calculation, the battery capacity parameter \(Q\) is expressed in terms of driving range (i.e., 340\,km) rather than energy units (kWh).  Fast-charging stations in the United States provide power outputs ranging from 50\,kW to 350\,kW. In this study, we use an average charging speed of 100\,kW for modeling purposes. Under this assumption, a full recharge of the 135\,kWh battery requires approximately 1.35 hours, equivalent to a replenishment rate of 251.9\,km per hour of charging. Thus, the charging rate parameter is set as \(\eta = 251.9\,\mathrm{km/hr}\). We also assume that trucks operate at a constant speed of 100\,km/h. The energy consumption rate \(\sigma\) is calculated as $ \sigma = 100\,\mathrm{km/hr} \times \frac{135\,\mathrm{kWh}}{340\,\mathrm{km}} = 39.7\,\mathrm{kWh/hr}$. The labor cost parameters are defined as follows: for leading trucks, the hourly cost is set to \(\alpha_1 = \$30/\mathrm{hr}\); for following trucks, we assign half of the leading truck cost, resulting in \(\alpha_2 = \$15/\mathrm{hr}\); and for dwell time, we use \(\alpha_3 = \$5/\mathrm{hr}\). To maximize the utilization of the leading truck rotation feature in this study, the additional cost for on-road platoon restructuring is set to 0.  For simplicity, all charging stations are assumed to have a uniform electricity price unless otherwise specified, just as in the small-scale test. According to a recent report \citep{chargeprice}, the average cost of fast charging at Level 3 or DC fast chargers in the U.S. ranges from \$0.40 to \$0.60 per kWh. We adopt the midpoint value of \$0.50/kWh for this study, and we initially set all charging stations with the same charging price, including destinations, unless otherwise noted, like in the small-scale test. Furthermore, the minimum and maximum allowable SOC are set to 0 and 1, respectively. Last, the maximum platoon size is set to \(L=4\). A summary of all parameter settings is provided in Table \ref{tab:para2}.

\begin{table}[htbp]
  \centering
  \caption{Parameter Setting}
  \label{tab:para2}
  \resizebox{\textwidth}{!}{
  \begin{tabular}{lll}
    \toprule
     Parameter & Description & Value\\
    \hline
    $L$ & Maximum platoon size & 4\\
    $\alpha_1$ & Unit labor cost respect to leading truck time & \$30/hr\\
    $\alpha_2$ & Unit labor cost respect to following truck time & \$15/hr\\
    $\alpha_3$ & Unit labor cost with respect to rest/dwell time at charging stations & \$5/hr\\
    $\alpha_4$ & Additional cost per platoon restructuring  & \$0 per swap\\
    $\gamma_i$ & The charging price at each charging station $i$ & $\$0.5/kWh$\\
    $\sigma$ & The energy consumption rate & $39.7 kWh/hr$\\
    $Q$ & Electric truck battery capacity & $340 km$  or $135 kWh$\\
    $\eta$ & Charging speed of charging stations & $251.9 km/hr$ or $100 kW$\\
    $SOC^u$ & Preset maximum allowable battery SOC after charging & 1\\
    $SOC^l$ & Preset lowest allowable SOC before charging & 0\\
    \bottomrule
\end{tabular}}
\end{table}

\subsection{Small-scale test}

\begin{figure}[htbp]
	\centering
    \includegraphics[scale=0.4]{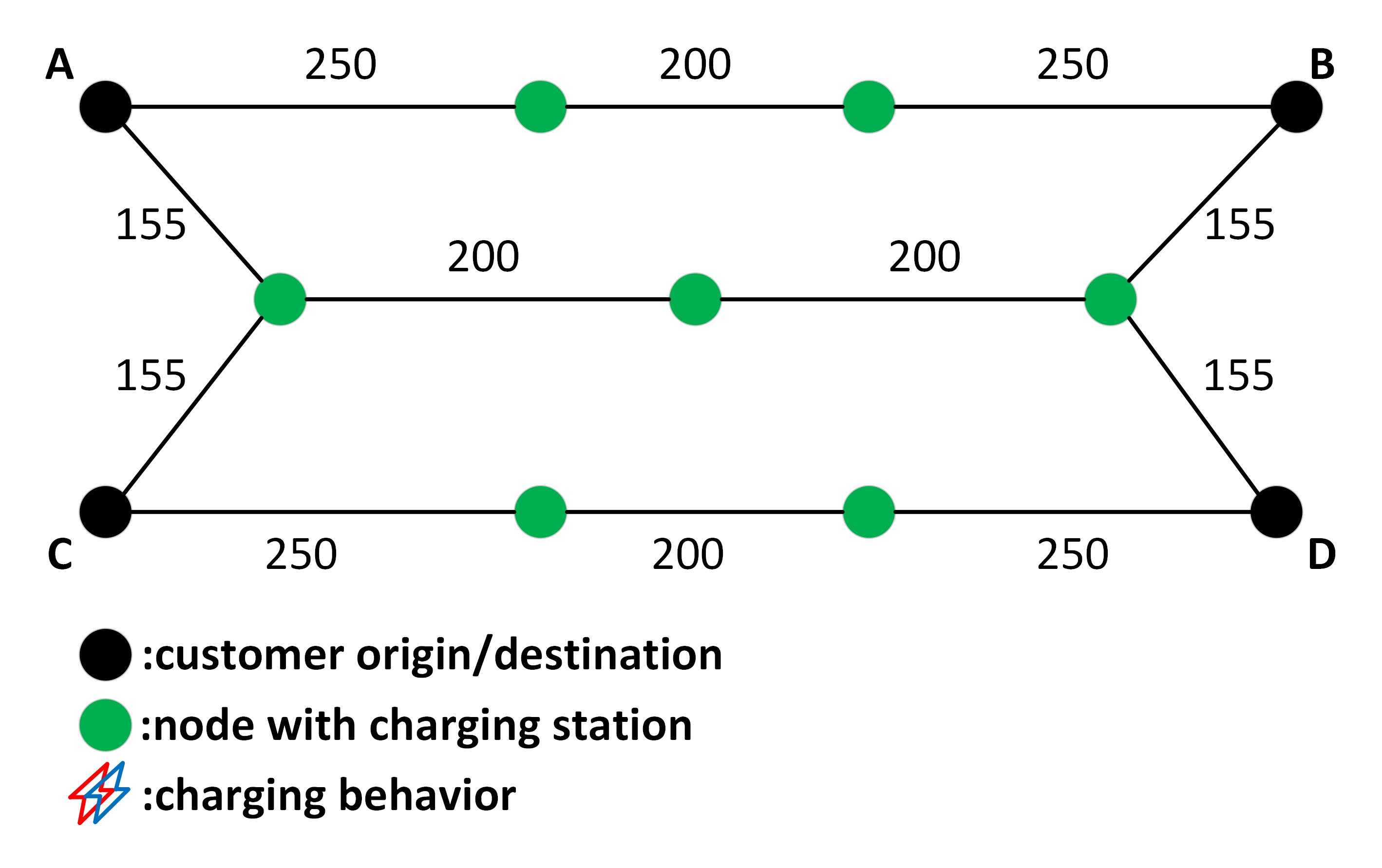}
	\caption{Small test network}
    \label{fig::evnetwork}
\end{figure}

In this subsection, we design a small test case with 2 truck deliveries on a network with 11 nodes, shown in Figure \ref{fig::evnetwork}. The network is modified based on the network shown in the illustrative example, Figure \ref{fig::platoonelectric}, by removing two charging stations from the second row. We use this small network because its refined size ensures us that an exact solution can be obtained by commercial solvers, such as CPLEX, in a timely manner. The test case is evaluated using the proposed mathematical model under two scenarios: one with platooning enabled and the other with platooning prohibited, representing traditional delivery operations. The latter scenario can be simply accomplished by our model when the platoon size limit is set to 1. We also introduce a new scenario when the platooning is enabled, but the leader swapping feature is prohibited, which means that trucks cannot change their positions in a platoon while traveling on road segments. Equivalently, it means that each truck's leading ratio, specified as $h_{i,j,k}$ is restricted to either 0 or 1, a binary variable, instead of any value between 0 and 1. Results comparison among these three scenarios will be presented in the remainder of this subsection to demonstrate the benefits of utilizing truck platooning technology in an electric truck delivery service.

To further explain the network shown in Figure \ref{fig::evnetwork}, we separate different kinds of nodes by distinct colors: black nodes are trucks' origins or destinations, and green nodes represent any intermediate nodes with charging stations. Specifically, the two trucks are departing from A and C and going to their corresponding destinations at B and D, respectively. To maximize the platoon potential and operation difference among the 3 scenarios, we set the latest arrival time of each destination node substantially large, [0,15], so that two trucks, even stops and charges at each charging station along the longest path, can still manage to meet the latest arrival time limits. In addition, to validate if the optimal solution matches our expectation on the platooning results shared in the illustrative example, we use platoon saving ratio as 0.15, the same as the illustrative example, and we set the charging price at the first and third charging station of the second row in Figure \ref{fig::evnetwork}  to be different from the others to reflect a higher charging cost at conjunctions, and their charging is set as the doubled price of the other charging stations. %Such action is also used to value the difference between charging once at the middle station and charging multiple times at all 3 stations, otherwise the charging cost is the same as we do not apply any additional cost for more stops. The difference matters because it will remark the benefit of leader swapping by making fewer stops.

Table \ref{evcomp} presents the operation cost for all three scenarios, along with the cost breakdown. We can witness that, when platooning is prohibited, referring to the traditional delivery service, the total operation cost is a bit higher than the scenario when platooning is allowed, no matter whether the leader swapping feature is enabled. In this test case, we decompose the operation cost into 3 categories, the charging cost, the on-road traveling cost that reflects the labor, and waiting cost that reflect the dwell time, and we can see that the traditional delivery service always yields the highest cost in each category. This tells that the traditional service usually needs more energy replenishment, and resulting in longer stops, and higher on-road traveling labor cost as there's no savings from less workload of followers. Specifically, the total savings yielded by platooning is $10.7\%$ $(=\dfrac{598.9-535}{598.9})$, and the charging-cost savings is $5.5\%$ ($=\dfrac{142.9-135}{142.9}$). In addition, comparing the platooning case without leader swapping, allowing the swapping of leading trucks can balance the energy consumption of different trucks and thus reduce the charging need and charging cost, which is $4.2\%$ ($=\dfrac{140.9-135}{140.9}$) in this example. This simple example highlights that platooning does provide substantial savings from different perspectives like energy consumption, labor cost reduction, and fewer or shorter en-route stops,  and we shall also never overlook the benefits from leader swapping feature as it acts as a further enhancement to the platooning technology, especially in energy re-balance in order to save charging costs and provide more flexibility to make stops at cheaper charging stations. Though the benefits amount may not be as large as the benefits showcased in this validation example, it is still worth deploying those operational strategies as the real-worth operational expense is considerable.

\begin{table}[htbp]
\caption{Comparison of operation costs with and without truck platooning}
\label{evcomp}
\centering

\resizebox{\columnwidth}{!}{\begin{tabular}{lccc}
\toprule

\textbf{Cost type} &  No-platoon & Platoon with leader swapping & Platoon w/o leader swapping    \\ 
\hline
\midrule
\textbf{Total operation cost} & 598.9 & 535.0 & 540.9  \\
\textbf{Traveling cost} & 420 & 366 & 366\\
\textbf{Waiting cost} & 36 & 34 & 34\\
\textbf{Charging cost} & 142.9 & 135 &  140.9\\
\hline

\bottomrule
\end{tabular}}
\end{table}

\begin{figure}[htbp]
	\centering
         \subfigure[Non-platooning case]{\label{fig:evnoplatoon}\includegraphics[scale=0.4]{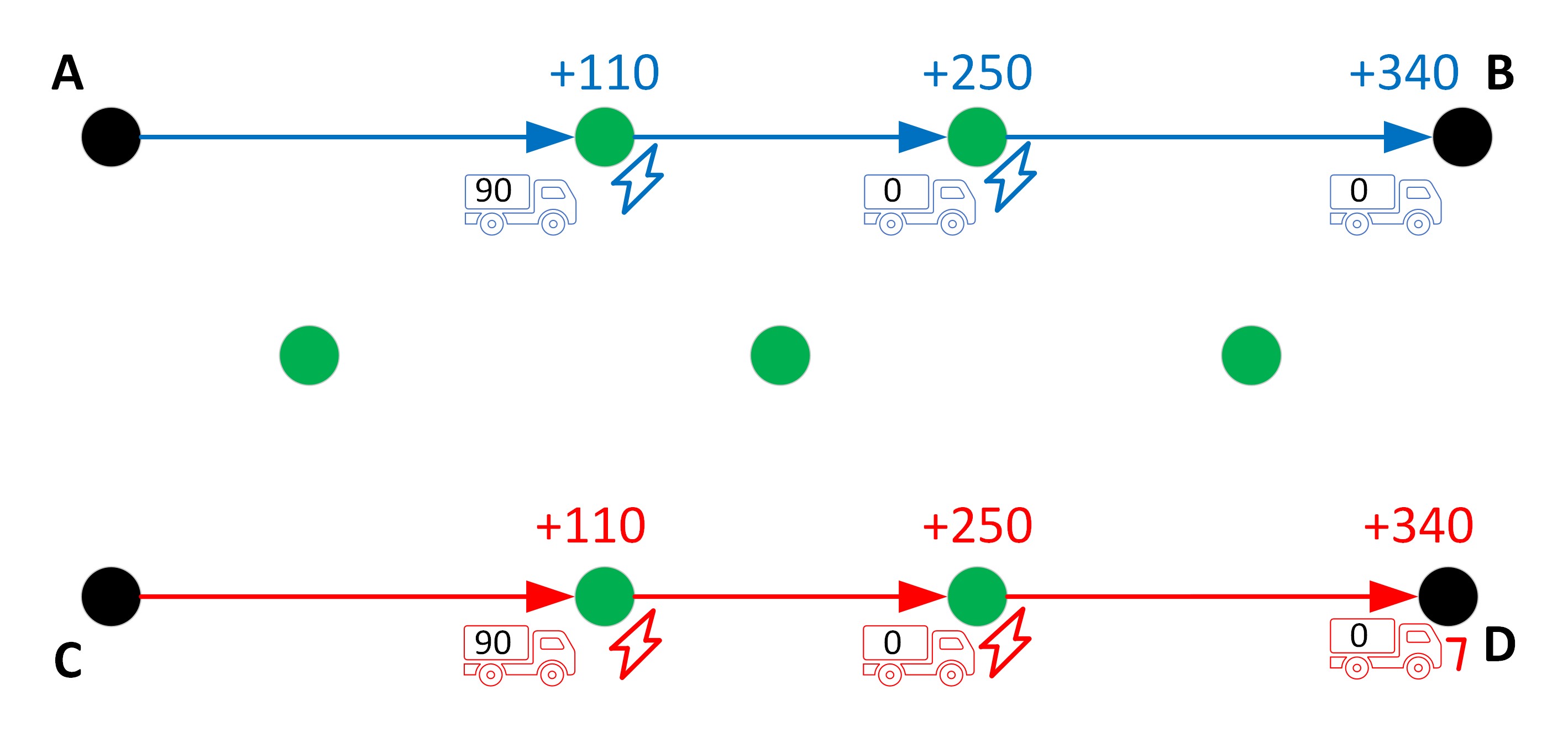}}\quad
	\subfigure[Platooning without leader swapping]{\label{fig:evplatoon}\includegraphics[scale=0.4]{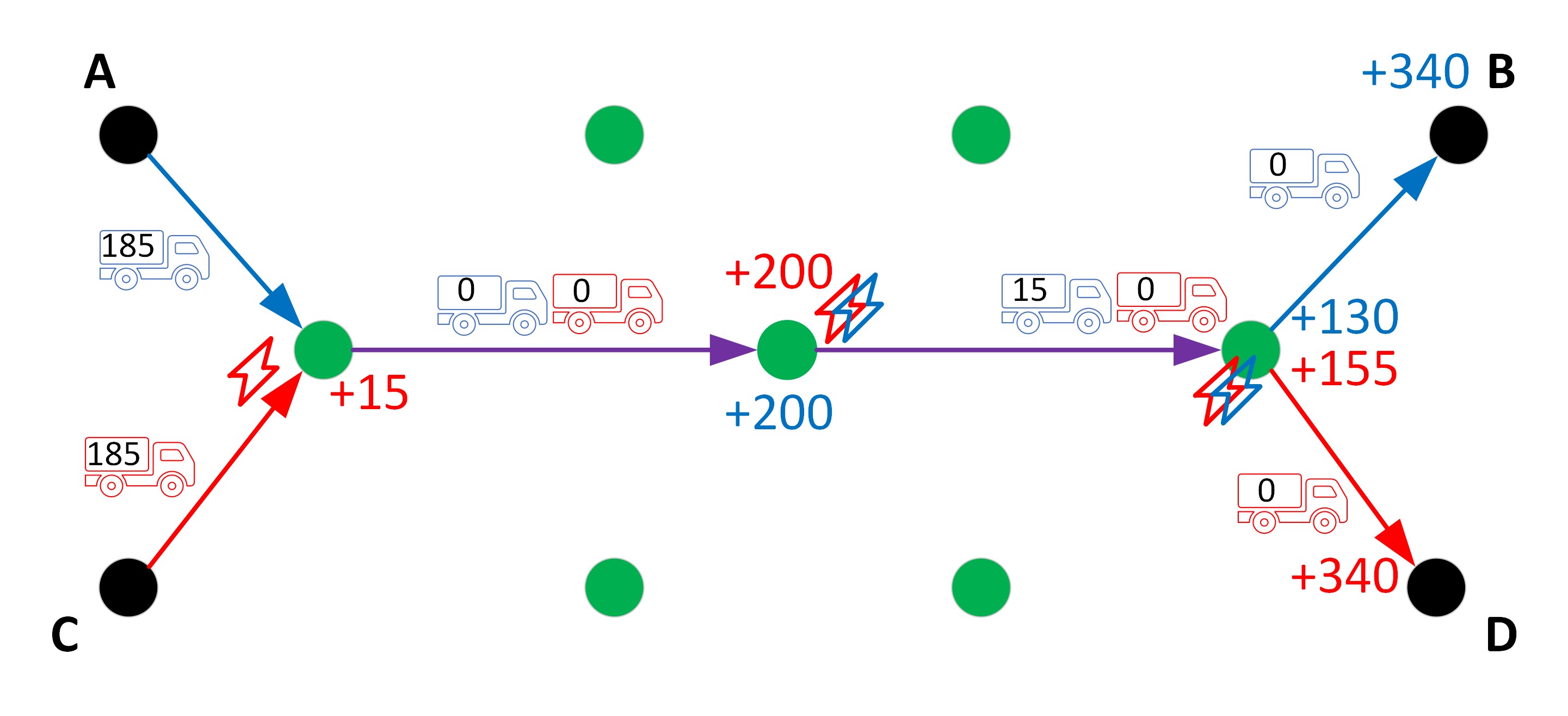}}\quad
    \subfigure[Platooning with leader swapping]{\label{fig:evplatoonswap}\includegraphics[scale=0.4]{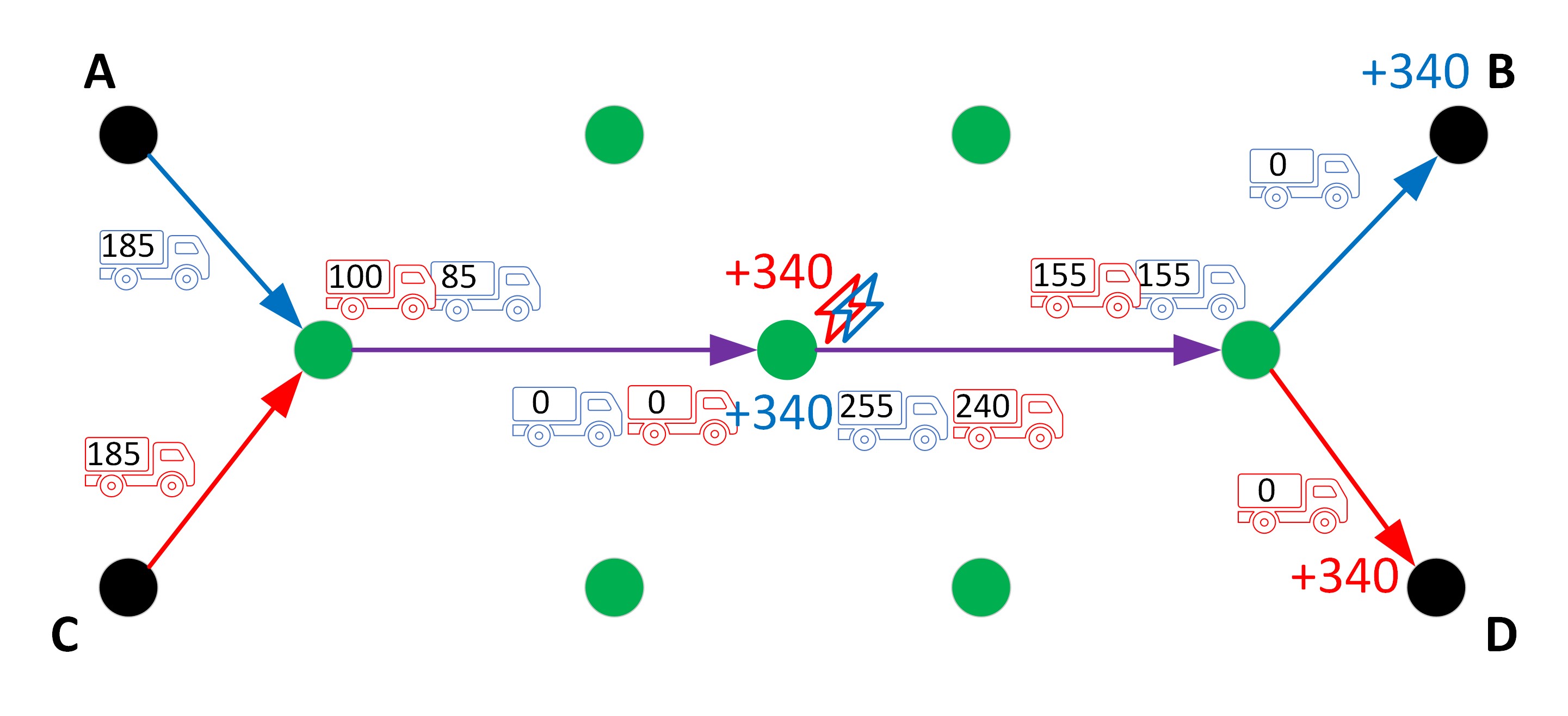}}\quad
    \caption{Operation comparisons of 3 different scenarios}
    \label{fig::evsmallresult}
\end{figure}

We further discuss the impact of platooning on operation plans, and visually present the operation plans under 3 scenarios in Figure \ref{fig::evsmallresult}. When platooning is not allowed, each truck will simply travel its shortest path to save the cost, and the resulting operation plan under this scenario is provided in Figure \ref{fig:evnoplatoon}, and the routes of different trucks are drawn with different colors and line formats. We further record the remaining driving range of each truck upon its arrival at next node along its traveling route, and the charging amount for each truck is also specified, with the lightning sign to indicate the charging behavior. For the second scenario when platooning is allowed but leader swapping is inactive, the operation plan is shown in Figure \ref{fig:evplatoon}. It can be seen that two trucks are now taking different routes and merge as a platoon on the common path. However, despite the symmetric network topology, the energy consumption between the two trucks is not evenly distributed. As a result, the truck departs from node C and arrives early at the first charging station to recharge the energy that supports 15km of driving range in order to have enough energy to reach the next charging station, because this truck needs to act as the leading vehicle and cannot benefit from the platooning. Similarly, two trucks have different charging amount at the last charging station that they visit. On the other hand, when the leader swapping is active, two trucks no longer need to make 3 stops along their routs. Instead, they can balance their energy consumption by alternatively acting as the leading vehicle and only charging once in the middle charging station with a full charge. The resulted operation plan is shown in Figure \ref{fig:evplatoonswap}, and swapping behavior can be found along the link with trucks showing up twice, representing that the two trucks swap their position in the platoon and therefore two different platoons are showcased. More specifically, each truck acts as the leader for 50\% of time while traveling along the common paths, and the numbers of 100 as well as 85 annotated in the Figure \ref{fig:evplatoonswap} presents the remaining mileage of two trucks upon reaching to the leader swapping moment. Based on the above small test example, we validate the optimal operation plans generated by our formulated MILP under 3 different scenarios, and the results emphasize that the platooning not only alter the route selection but also the charging decisions of each truck. The influence of the leader swapping feature is further highlighted by the different charging amount and charging stops of the two trucks, which brings additional improvement to the existing platooning technology under mid-haul electric truck logistics. 

\subsection{Large-scale test}

\begin{figure}[htbp]
	\centering
    \includegraphics[width=\textwidth]{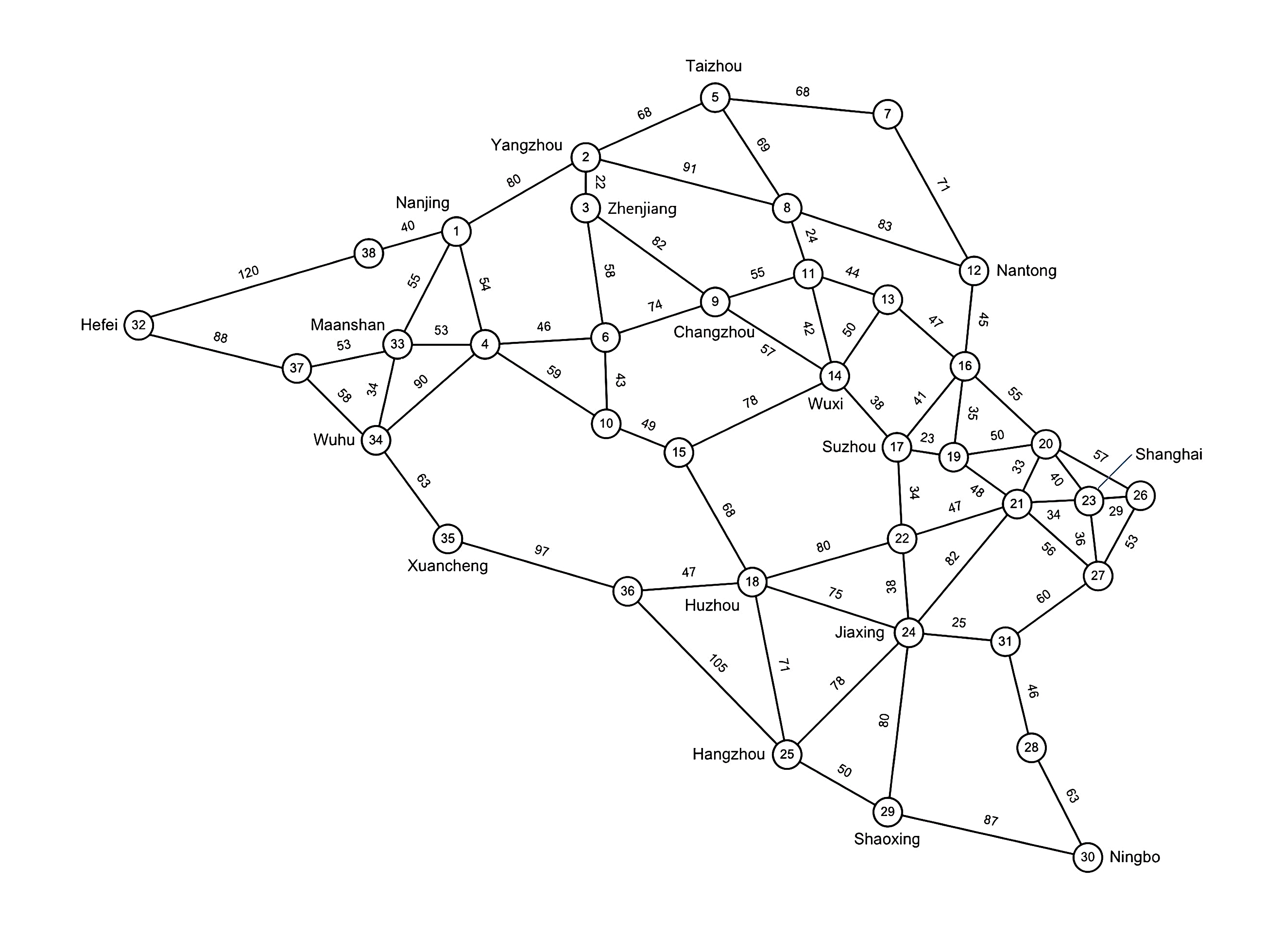}
	\caption{Yangtze River Delta network \citep{delta2022}}
    \label{fig::yang2}
\end{figure}

In this subsection, we provide a numerical experiment on a real-world transportation network, the Yangtze River Delta network. It is located in the central eastern area of China, an economic center where logistics play an essential role. The network consists of 38 nodes or cities, and the whole network is presented in Figure \ref{fig::yang2} with distances between nodes labeled. The units of labeled distances are in kilometers, so we manually transform them into time dividing by an average truck speed on highways, which is 100km per hour as also mentioned before. In addition, the platoon saving ratio, $\beta$, is set as 0.1.

For the test case design, we still set the identical truck battery capacity as 135kWh which allows each truck to travel 340km at most. Trucks' OD pairs are randomly generated from the network with preset random seeds (seed number as 1234 plus the execution number with Java SE 8 for reproducibility), so are their latest arrival time, and their latest arrival time is randomly chosen between the shortest travel time from origin to destination and the largest value of the shortest travel times associated with every possible OD pairs with a uniform distribution. For the first experiment to compare CPLEX and our proposed heuristic solution algorithm, we generate 5 different cases for each specific problem size, with each case assigned a fixed random seed for regeneration. In addition, we set a 3600-second running time limit for CPLEX to output a feasible solution while not wasting too much time. The time limit for the heuristic is also set as 3600 seconds, or no improvement in the solution in 50 consecutive iterations. The running time and result comparisons between CPLEX and our proposed algorithm are presented in Table \ref{EVYRDresult}.

\begin{table}[htbp]
\caption{Comparative results with CPLEX}
\label{EVYRDresult}
\centering
\begin{threeparttable}
\resizebox{\columnwidth}{!}{\begin{tabular}{lllllll}
\toprule
\multicolumn{2}{l}{\textbf{Instance}}& \multicolumn{5}{l}{\textbf{Result}} \\
\cmidrule(lr){1-2} \cmidrule(lr){3-7}
\makecell{\textbf{Number of} \\ \textbf{trucks}} & \textbf{Execution count}  & \makecell{\textbf{CPLEX solution}} & \makecell{\textbf{CPLEX's} \\ \textbf{CPU time}}  & \textbf{Heuristic solution} & \makecell{\textbf{Heuristics's} \\ \textbf{CPU time}}   &\textbf{Gap} \\ 
\hline
\midrule
5 & 1 & 291 & 104.9s & 296 & 1.2s & 1.72\% \\
& 2& 600 & 3,600s & 602& 8.5s & 0.33\% \\
& 3& 457 & 3,600s & 461 & 0.9s & 0.87\% \\
& 4& 538 & 3,600s & 540 & 3.8s & 0.37\% \\
& 5& 583 & 3,600s & 581 & 7.7s & -0.34\% \\
6 & 1& 400 & 540.8s & 415 & 0.8s & 3.75\%\\
& 2& 638& 3,600s & 698 & 7.0s & 9.40\% \\
& 3& 490 & 3,600s & 527 & 1.3s & 7.55\% \\
& 4& 492& 3,600s & 493 & 3.9s & 0.20\% \\
& 5& 588& 3,600s & 605 & 6.0s & 1.19\% \\
7 & 1& 464 & 3,600s & 481 & 1.3s & 3.66\%\\
& 2& 700& 3,600s & 712 & 2.1s & 1.71\% \\
& 3& 614& 3,600s & 626 & 0.9s & 1.95\% \\
& 4& 567& 3,600s & 595 & 4.7s & 4.94\% \\
& 5& 606 & 3,600s & 607 & 7.1s & 0.17\% \\
8 & 1& 562 & 3,600s & 610 & 1.6s & 8.54\% \\
& 2& 833& 3,600s & 890 & 2.1s & 6.84\% \\
& 3& 652 & 3,600s & 675 & 3.4s & 3.53\% \\
& 4& 658 & 3,600s & 690 & 5.5s & 4.86\% \\
& 5& 650 & 3,600s & 623 & 7.3s & -4.15\% \\
9 & 1& 590 & 3,600s & 604 & 1.7s & 2.37\% \\
& 2& N/A & 3,600s & 1038 & 2.4s & N/A \\
& 3& 684 & 3,600s & 669 & 1.8s & -2.19\% \\
& 4& 700 & 3,600s & 740 & 4.6s & 5.71\%\\
& 5& 766 & 3,600s & 763 & 1.9s & 0.39\% \\
10 & 1& 695 & 3,600s & 671 & 2.1s & -3.45\%\\
& 2& N/A & 3,600s & 1152 & 2.8s & N/A \\
& 3& 738& 3,600s & 746 & 3.4s & 1.08\% \\
& 4& 836& 3,600s & 814 & 4.0s & -2.63\% \\
& 5& 824& 3,600s & 836 & 2.2s & 1.46\% \\
\hline

\bottomrule
\end{tabular}}
\begin{tablenotes}
    \footnotesize
    \item Gap: the relative gap between CPLEX and heuristic results,$\dfrac{Obj_{heuristic}-Obj_{CPLEX}}{Obj_{CPLEX}}$.
\end{tablenotes}
\end{threeparttable}
\end{table}

According to Table \ref{EVYRDresult}, CPLEX successfully identified the optimal solution within the 3600-second time limit in only 2 out of 25 executions, specifically, in execution 1 of the 5-truck case and execution 1 of the 6-truck case. While CPLEX was able to generate high-quality solutions for most instances within the time threshold, its ability to guarantee optimality diminishes as problem size increases. In certain larger cases, such as the second executions of the 9-truck and 10-truck instances, CPLEX fails to produce even a feasible solution within 3600 seconds, as indicated by the “N/A” entries in Table \ref{EVYRDresult}.

In terms of solution quality, CPLEX generally achieved slightly better results than ALNS, with an average relative gap of 2.6\%, computed as the mean of all 23 valid gaps presented in the last column of Table \ref{EVYRDresult}. However, the performance difference is modest. In addition, although CPLEX outperforms ALNS in small instances (e.g., 5- and 6-truck cases), its advantage shrinks as the problem size grows. In fact, in 5 out of the last 10 executions for larger instances, ALNS produces better solutions than CPLEX, highlighting CPLEX’s scalability limitations in handling this complex problem. By contrast, the proposed ALNS method demonstrates clear superiority in computational efficiency. For small cases, ALNS typically finds solutions within a few seconds, reducing CPU time by over 99\% compared to CPLEX—while still delivering solutions with relative gaps often less than 5\%. This efficiency becomes increasingly critical in larger-scale real-world scenarios where the solution space expands exponentially with the number of trucks and network size. Overall, these results underscore the strong scalability and practical utility of the proposed ALNS approach in generating high-quality solutions within significantly shorter runtimes.

To further evaluate the scalability of the proposed ALNS algorithm in handling larger instances and better reflecting real-world logistics scenarios, we set up a series of test instances with fleet sizes ranging from 5 to 150 trucks. In these experiments, truck origins, destinations, and time windows are randomly generated. The time windows are set large enough to ensure the feasibility of each instance. The results of these scaling tests are presented in Table \ref{ScaleresultET}. Similarly, the time limit in this experiment is set as 3600 seconds, and the random cases are also solved using CPLEX to evaluate the solution quality of the ALNS algorithm. 

In addition, we assess the platooning benefit observed in each randomly generated case. Specifically, the platooning benefit is quantified as the difference in total operation costs between scenarios with and without platooning. To compute this, each instance is run twice: once with a platoon size limit of 4, and once with a platoon size limit of 1, representing the non-platooning scenario. Both instances use the results obtained through ALNS as CPLEX may not be able to yield a feasible solution for large instances.  The cost difference between these two runs reflects the aggregate benefit of platooning, which in this context arises not only from reductions in energy consumption but also from savings in on-road labor costs due to decreased driver workload. Moreover, forming larger and more platoons can potentially reduce dwell times by lowering the charging time.

The results in Table \ref{ScaleresultET} demonstrate that ALNS ensures solution quality when compared to CPLEX, especially when CPLEX manages to output good solutions, as shown in the first three rows of  Table \ref{ScaleresultET}. In the first row, when CPLEX solves to optimality with 3161 seconds, ALNS only underperforms by 0.86\%. In addition, ALNS outperforms CPLEX by 0.14\% and 8.81\% respectively for the 10-truck and 15-truck cases. For larger cases, CPLEX either fails to provide a feasible solution within the time limit or runs into the ``Out of memory" error. We tried our best to test CPLEX on different setups, including WorkMem, NodeFileInd, and MIPGap in order to avoid such errors, but the test instance is still too complicated for CPLEX to solve. 

In terms of computational efficiency, ALNS exhibits strong scalability as the problem size increases. Even for the largest instance with 150 trucks, the algorithm achieves convergence in approximately 120 seconds. Even for the small case with only 5 trucks, ALNS saves $\dfrac{3161 - 1.2}{3161} = 99.96\%$ of computation time while still providing a competitive solution. This level of computational efficiency supports practical applications where multiple runs under varying parameter settings may be necessary. Regarding platooning benefits, the results indicate an upward trend as the problem size grows because more trucks offer more platoon potential, leading to more significant savings. Notably, the platooning benefit reaches up to 2.77\% in the 110-truck case.

\begin{table}[htbp]
\caption{Results of scaling instances}
\label{ScaleresultET}
\centering
\begin{threeparttable}
\resizebox{\columnwidth}{!}{\begin{tabular}{llllllll}
\toprule
\multicolumn{1}{l}{\textbf{Instance}}& \multicolumn{4}{l}{\textbf{Result}} \\
\cmidrule(lr){1-1} \cmidrule(lr){2-8}
\textbf{Number of trucks} & \textbf{CPLEX time}  & \textbf{CPLEX solution} & \textbf{Heuristic time}  & \textbf{Heuristic solution} &\textbf{Gap} &\textbf{Platoon benefit} & \textbf{\% Platoon benefit}\\ 
\hline
\midrule
5& 3,161s & 465 & 1.2s & 469 &0.86\% & 0 & 0\%\\
10& 3,600s& 712 & 4.3s & 711 & -0.14\% & 0 & 0\%\\
15& 3,600s & 1,748 & 4.3s & 1,594 & -8.81\% & 19& 1.18\%\\
20& 3,600s & N/A & 7.2s & 2,177 &  N/A &0 & 0\%\\
25& 3,600s & N/A & 11.3s & 1,629 & N/A & 77.34& 4.53\%\\
30& 3,600s & N/A & 19.3s & 4,138 & N/A & 90 & 2.13\%\\
35& 3,600s &  N/A & 16.9s & 2,083 & N/A & 42 & 1.98\%\\
40& Out of memory & N/A & 13.2s  & 5,693 &  N/A& 101 & 1.74\% \\
50& Out of memory& N/A & 15.8s  & 7,051 &  N/A& 79 & 1.11\%\\
60& Out of memory& N/A &  46.1s & 8,107 &  N/A & 132 & 1.60\%\\
70& Out of memory& N/A &  58.1s & 9,947 &  N/A & 199 & 1.96\%\\
80& Out of memory& N/A & 111.9s &  10,543 &  N/A & 117 & 1.10\%\\
90& Out of memory& N/A &  76.3s & 12,836 &  N/A & 204 & 1.56\%\\
100& Out of memory& N/A &  80.3s &  13,715 &  N/A & 268 & 1.92\%\\
110& Out of memory& N/A &  91.2s &  11,143 &  N/A & 317 & 2.77\%\\
120& Out of memory& N/A &  91.3s &  11,107 &  N/A & 252 & 2.22\%\\
130& Out of memory& N/A &  97.1s &  12,087 &  N/A & 301 & 2.43\%\\
140& Out of memory& N/A &  109.1s &  13,654 &  N/A & 212 & 1.53\%\\
150& Out of memory& N/A &  118.7s &  15,079 & N/A & 384 & 2.48\%\\
\hline

\bottomrule
\end{tabular}}
\begin{tablenotes}
    \footnotesize
    \item Gap: the relative gap between CPLEX and heuristic results, $\dfrac{Obj_{heuristic}-Obj_{CPLEX}}{Obj_{CPLEX}}$.
    \item Platoon benefit: the difference between no-platoon and platoon results by heuristic
    \item \% Platoon benefit: the platoon benefit in percentage, $\dfrac{Platoon\ benefit}{Objective\ without\ platoon}$.
    \item N/A: No solution found within the time limit.
\end{tablenotes}
\end{threeparttable}
\end{table}

\subsection{Sensitivity Analysis}
In this subsection, sensitivity analyses are conducted to quantitatively represent the impact of key parameters on the results. The key parameters 
 are the platoon saving ratio, truck battery capacity, and labor cost coefficient in leading trucks.

 For all the experiments in the sensitivity analysis section, it is noted that they are tested on a 4 by 4 grid network with an identical edge length of 150km instead of the complicated Yangtze River Delta network so that an optimal solution can be secured by CPLEX within a reasonable timeframe like 10800s. In addition, this choice is primarily motivated by the inherent stochasticity of the ALNS algorithm compared to decomposition-based approaches. In large-scale real-world networks, if the parameter variations are not sufficiently significant: for example, increasing the platoon saving ratio from 0.05 to 0.075, then the impact of such minor changes may be obscured by the randomness of the heuristic search process, which can be reflected by the operator selection mechanism and solution acceptance criteria that are highly stochastic. To avoid such bias, the sensitivity analysis is therefore conducted on a smaller, controlled network where the problem can be solved to optimality, allowing the influence of parameter adjustments to be observed more clearly. Additionally, we chose the network size as 4 by 4 because according to our preliminary tests, CPLEX may encounter the ``out of memory" issue or run time error for a 5 by 5 grid network. 
 
\subsubsection{Platoon saving ratio}
\begin{figure}[htbp]
	\centering
	\subfigure[Operation cost]{\label{fig:evbetaa}\includegraphics[scale=1]{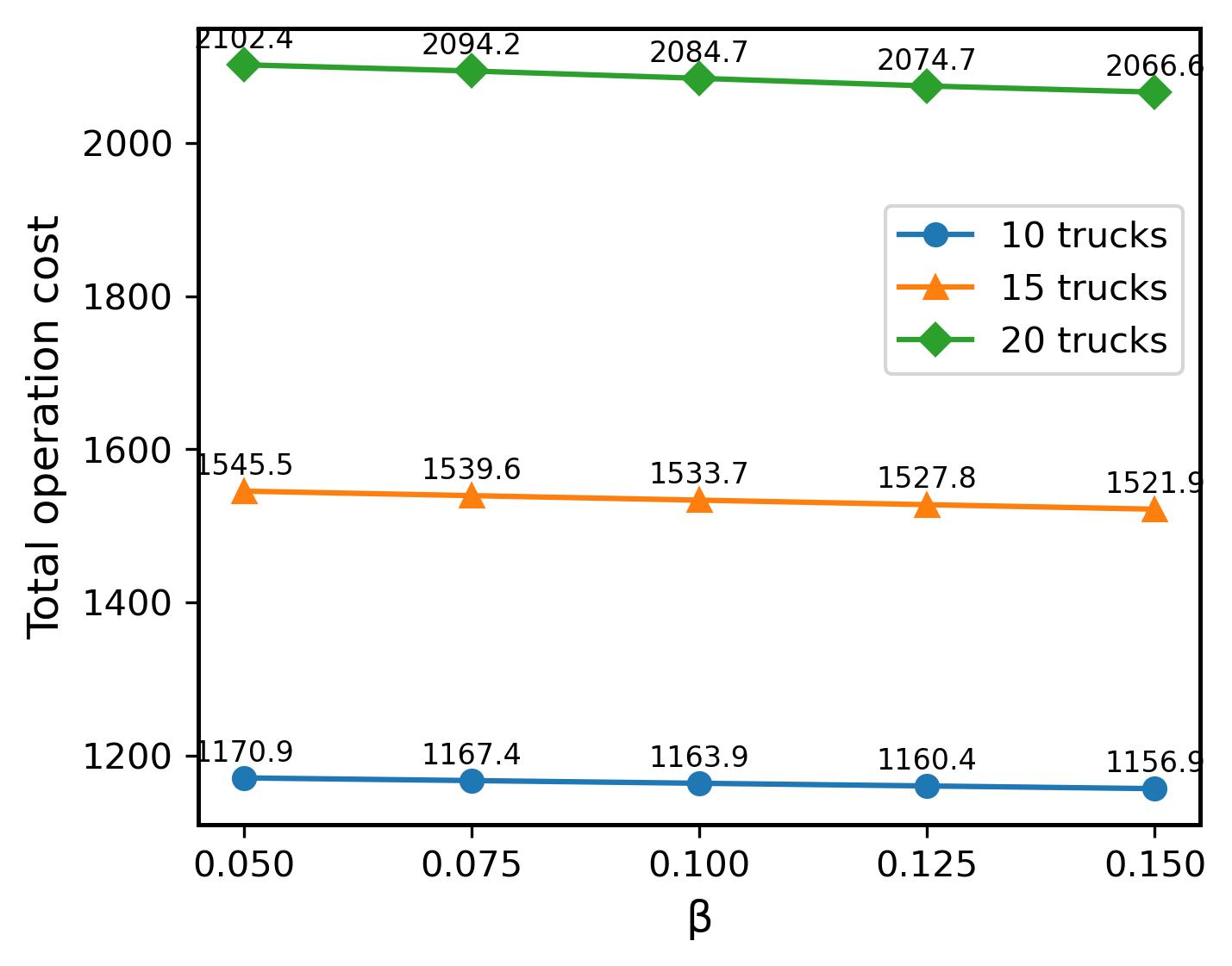}}\quad
	\subfigure[Cost breakdown of 20-truck case]{\label{fig:evbetab}\includegraphics[scale=1]{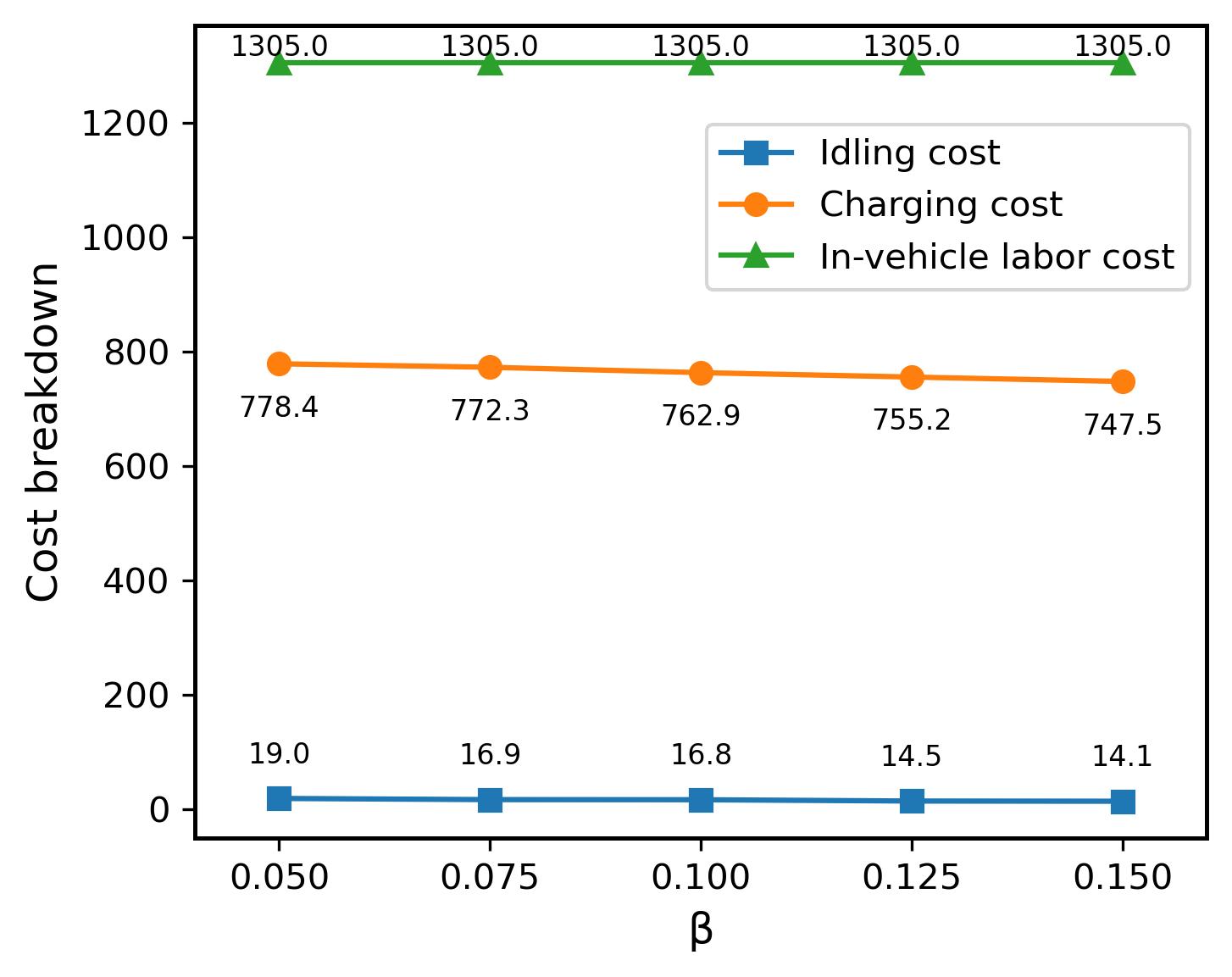}}\quad\\
	\caption{Influence of platoon saving ratio on all types of costs}
    \label{fig::evbeta}
\end{figure}

The platoon saving ratio is often selected to be 0.1 or 0.05 for the following trucks in a platoon, but with the advancement of better aerodynamic design or autonomous driving to reduce the headways, such a ratio is anticipated to be larger. Therefore, Figure \ref{fig::evbeta} plots the operation cost with different platoon saving ratios under different delivery demands. Specifically, in Figure \ref{fig:evbetaa}, we plot the overall trend of the total operation cost for 3 different cases: 10, 15, and 20 trucks. As expected, no matter under what level of delivery demand, an increase in the saving ratio will always yield a reduction in the operation cost, as more energy can be saved in the same operation plan. It can also be seen that the relationship between the operation cost and the platoon saving ratio is exactly linear for all cases. In addition, as we set the increment for the platoon saving ratio to be identical as 0.025, the identical increment in the operation cost infers that the platoon formations in our cases have not been altered within the changes of platoon saving ratio. One possible explanation is that the limited increase from 0.1 to 0.15 or decrease from 0.1 to 0.05 does not provide any alternative paths that can save more compared to the current ones. 

Diving into the cost breakdown shown in Figure \ref{fig:evbetab}. We notice that the in-vehicle labor cost, that is associated with leading and following travel time, does not change at all, which tells us that the routing decisions are not impacted by $\beta$ changes. However, as the major benefit of platooning, the energy consumption decreases, which in turn lowers the charging cost. Idling cost with respect to dwell time, in the meanwhile, also shrinks because of less charging time en route.

\subsubsection{Truck battery capacity}

As the primary motivation for investigating the impact of platooning on electric trucks arises from their range anxiety, it is essential to explore how these impacts vary with changes in battery capacity. This analysis is particularly meaningful given the expectation that future electric trucks will increasingly be equipped with larger batteries to support extended driving ranges.

\begin{figure}[htbp]
	\centering
    \includegraphics[scale = 1]{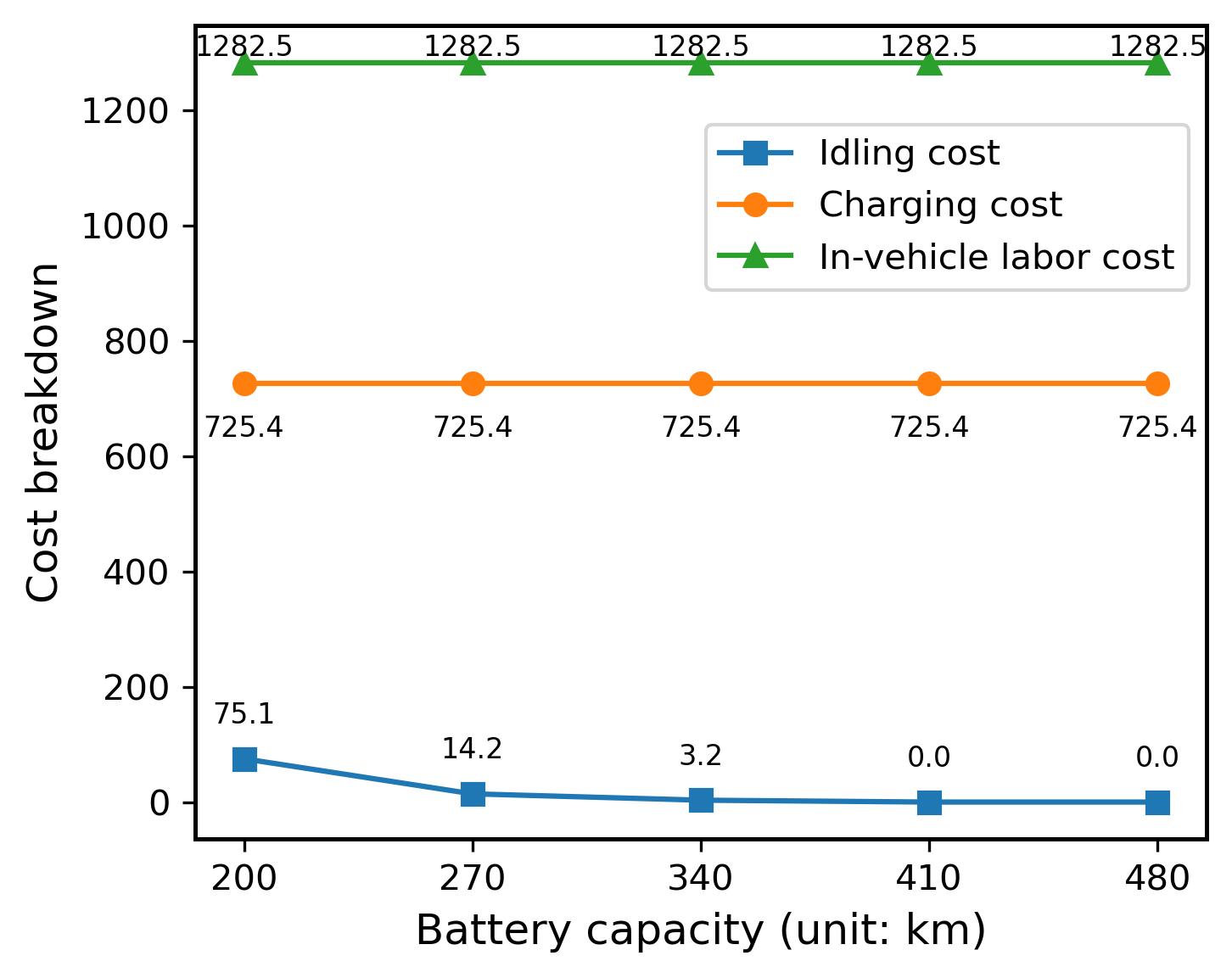}
	\caption{Influence of truck capacity on all types of costs}
    \label{fig::evQ}
\end{figure}

Starting from the baseline battery capacity of 135 kWh, which enables a driving range of 340 km, we examine four alternative capacities: 79.4kWh, 107.2 kWh, 162.8 kWh, and 190.6 kWh. These correspond to driving ranges of 200 km, 270 km, 410 km, and 480 km, respectively. The increments in driving range are set to be uniform to facilitate a clearer analysis of the trend. Additionally, in Figure \ref{fig::evQ}, the results are presented in terms of driving range (km) rather than energy capacity (\,kWh) to enhance the interpretability of the findings. To be noted, the cost breakdown plot is based on the 20-truck case. 

It is evident from Figure~\ref{fig::evQ} that the in-vehicle labor cost remains constant across all battery capacity scenarios, indicating that the routing plans remain unchanged in terms of optimality. Similarly, the charging cost is consistent across different battery capacities. This is because, with all other parameters held constant, such as the platoon saving ratio and energy consumption rate, traversing the same routes leads to identical total energy consumption, regardless of battery size.

In contrast, a significant reduction in waiting cost is observed, dropping from 75.1 to 0 as the battery-supported driving range increases from 200 km to 410 km. This reduction occurs because smaller battery capacities necessitate more frequent charging stops, resulting in increased dwell times and thus higher waiting costs. This phenomenon highlights the operational advantage of equipping electric trucks with larger batteries, as it minimizes charging interruptions during trips. Furthermore, it is notable that when the battery capacity supports driving ranges beyond 410 km, the waiting cost reduces to zero. This implies that trucks no longer require intermediate charging en route and can complete their trips with only destination charging, thereby eliminating any charging-induced delays.

\subsubsection{Unit labor cost of leading trucks}
\begin{figure}[htbp]
	\centering
    \includegraphics[scale = 1]{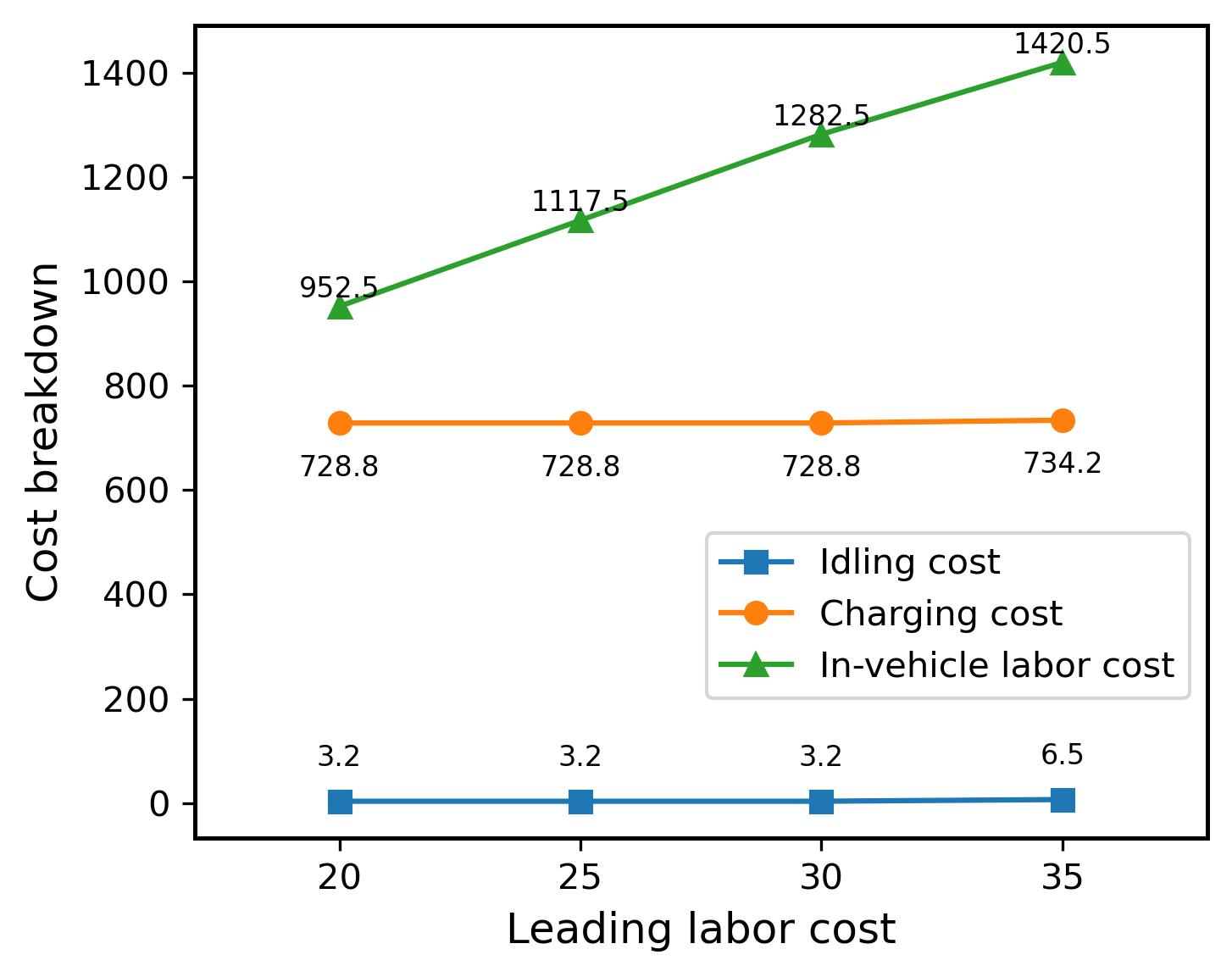}
	\caption{Influence of leading labor cost on all types of costs}
    \label{fig::evalpha1}
\end{figure}

Labor costs are a critical component of logistics operations, and it is therefore essential to investigate how variations in labor rates influence overall costs. This is particularly relevant given the substantial differences in wage levels across countries and regions. In this subsection, we examine four distinct unit labor cost levels for drivers operating leading trucks, as this reflects the full workload borne by these drivers. The experiments are conducted on a 20-truck scenario, and the resulting cost breakdowns at each unit labor cost level, \$20, \$25, \$30, and \$35 per hour, are visually presented in Figure \ref{fig::evalpha1}.

According to Figure \ref{fig::evalpha1}, the charging cost remains unchanged across the first three unit labor cost levels (\(\alpha_1 = 20, 25, 30\)). This suggests that the charging decisions are not significantly influenced by variations in the leading truck labor cost within this range. One contributing factor is that the labor cost associated with dwell time (\(\alpha_3\)) remains constant throughout the experiments. Additionally, the travel cost exhibits a consistent increase of 165 units when \(\alpha_1\) rises from 20 to 25 and from 25 to 30. However, this increment decreases to 138 units when \(\alpha_1\) increases from 30 to 35. This reduction in the rate of increase likely indicates that, beyond a certain labor cost threshold, the system begins to adjust routing strategies to mitigate the higher expenses associated with leading truck operation. Specifically, trucks may be rerouted along alternative paths with shorter travel times or form more or larger platoons. These adjustments reduce the exposure of individual trucks to leading roles or solo travel, thereby lowering labor-related travel costs despite rising wage levels. Furthermore, both the waiting cost and the charging cost show noticeable changes at \(\alpha_1 = 35\) compared to the earlier levels. This indicates that modifications in routing and platoon formation simultaneously influence charging decisions and dwell times. Such interdependencies highlight the complex trade-offs between labor costs, routing plan, and energy management in electric truck platooning operations.

\section{Conclusion and future direction}\label{sec:con}
In this study, we investigated the integration of truck platooning into electric truck logistics with charging consideration. This study highlights the significance of platooning as a strategy to alleviate range anxiety through energy savings and how the energy savings brought by platooning can be used to enhance electric trucks' routing and charging decisions by providing alternative paths that capitalize on platooning benefits. Unlike previous studies on electric truck platooning, which often simplify the problem to a single highway corridor, our model extends the scope to an expanded road network to more accurately represent real-world conditions. Additionally, we introduce an innovative leader-swapping mechanism within platoons to balance energy consumption among trucks during shared travel segments. The proposed MILP model jointly optimizes routing, charging, and platooning decisions to minimize total operation costs, including charging expenses and labor costs, while accommodating varying time windows and heterogeneous charging prices.

To effectively solve this complex problem, particularly for large-scale instances, we developed a solution framework combining a customized pre-processing procedure with an Adaptive Large Neighborhood Search (ALNS) algorithm. The pre-processing phase identifies candidate charging stations for each truck to minimize operation costs and constructs a high-quality initial solution, which serves as a warm start for ALNS. Within ALNS, the problem is reformulated into encoded truck segments that encapsulate routing, platooning, and charging information. A diverse set of removal and insertion operators is employed to explore routing and platooning possibilities, while swapping operators, such as leading ratio adjustments and charging amount modifications, optimize charging decisions and adjust platoon formations. Furthermore, an adaptive weight update mechanism and a simulated annealing-based solution acceptance criterion are incorporated to iteratively improve solution quality and avoid premature convergence to local optima. This framework demonstrates the ability to generate high-quality solutions for large-scale problems within a computationally efficient timeframe.

We first conduct a small-scale experiment to validate our mathematical model and to illustrate the impact of platooning and leader swapping on electric truck operations. Subsequent large-scale experiments on the Yangtze River Delta network are performed to assess the performance and scalability of the heuristic algorithm. Results indicate that the ALNS method can deliver high-quality solutions while achieving up to 99.96\% savings in computation time compared to CPLEX. For instances with more than 30 trucks, CPLEX fails to find feasible solutions within the time limit, whereas ALNS successfully solves these cases within 120 seconds, demonstrating its scalability and practical utility in real-world applications. Moreover, scalability tests reveal notable platooning benefits in large instances, with cost savings reaching up to 2.77\%. Sensitivity analyses are also performed to assess the effects of key parameters on operation costs. The results show that the platooning saving ratio significantly influences routing plans but has little impact on charging decisions. In contrast, battery capacity strongly affects charging strategies, and unit labor costs associated with leading truck operation can influence routing decisions, triggering cascading effects on charging schedules.

To further improve our research on electric truck platooning, there are several potential directions. From the perspective of problem formulation, incorporating more comprehensive time window constraints, such as earliest and latest departure times, would better accommodate time-sensitive operational needs. In addition, several simplifying assumptions made in this study could be relaxed to more accurately reflect real-world scenarios. For example, modeling variable truck battery capacities, heterogeneous energy consumption rates, and non-linear or station-specific charging speeds could provide a richer and more realistic depiction of electric truck operations. These factors merit further investigation, as they have the potential to significantly influence routing, scheduling, and charging decisions, and may amplify or diminish the overall benefits of platooning. In terms of the problem context, future work can also explore integrating additional elements relevant to electric truck operations. These include the deployment of battery swapping facilities to accelerate energy restoration, modeling charging queues under high-demand scenarios, and incorporating mandatory rest breaks for drivers, which are critical in mid-haul logistics. Exploring how these factors interact within a unified framework could yield valuable insights, though it would inevitably increase the complexity of the problem. From the solution methodology perspective, one promising direction for enhancing the performance of the current ALNS approach lies in systematic hyperparameter tuning. Similar to the practices used in machine learning, it would be beneficial to collect a diverse set of benchmark instances with known optimal solutions and apply techniques such as grid search or Bayesian optimization to identify the best settings for key ALNS parameters. These include the simulated annealing temperature schedule, the $\delta$ values in the adaptive weight update procedure, and convergence criteria such as the maximum number of iterations without improvement. With more comprehensive hyperparameter tuning, it is anticipated that ALNS could achieve faster convergence while improving solution quality.

\bibliographystyle{apalike} 
\footnotesize\bibliography{mybib}
\end{document}